\documentclass[oneside]{article}
\usepackage[utf8]{inputenc}
\usepackage[english]{babel}  

\usepackage{amsthm}

\theoremstyle{definition}

\usepackage{amsmath}
\usepackage{amsfonts}
\usepackage{amssymb}
\usepackage{graphicx}
\numberwithin{equation}{section}
\usepackage{lipsum}
\usepackage{mathtools}
\usepackage{scalefnt}
\usepackage{empheq}
\usepackage[a4paper,left=3.5cm,right=3cm,top=3cm,bottom=3cm]{geometry}
\usepackage{makeidx}
\usepackage{bm}
\makeindex

\usepackage{indentfirst}

\usepackage{fullpage}
\usepackage{amsmath}
\usepackage{amssymb}
\usepackage{graphicx}
\usepackage{subfigure}
\usepackage{color}
\usepackage{epsfig}
\usepackage{epstopdf}
\usepackage{hyperref}
\usepackage{algorithm,algcompatible}
\usepackage[affil-it]{authblk}
\usepackage{comment}

\date{}

\author[1]{Thi Tam Dang\thanks{\texttt{tam.dang@helsinki.fi}}}
\author[2]{Trung Hau Hoang\thanks{Corresponding author: \texttt{trunghaugg@gmail.com}}}

\affil[1]{Department of Mathematics and Statistics, University of Helsinki, Finland}
\affil[2]{Department of Mathematics, University of Innsbruck, Austria}

\begin{document}
	
	\title{From Memory Model to CPU Time: Exponential Integrators for Advection-Dominated Problems} 

	\maketitle
	
\begin{abstract}
In this paper, we investigate the application of exponential integrators to advection-dominated problems. We focus on Krylov subspace and Leja interpolation methods to compute the action of exponential and related matrix
functions. Complementing our earlier paper, “Should Exponential Integrators be Used for Advection-Dominated Problems?” (to appear in Advances in Applied Mathematics and Mechanics, 2025) based on a performance model, we extend the numerical investigation to higher-order Krylov approximations and new numerical regime, and assess their CPU-time efficiency relative to explicit Runge–Kutta schemes. We show that, depending on the problem setting, exponential integrators can either outperform or match explicit Runge–Kutta schemes. We also observe that Leja-based methods outperform Krylov iterations for large time steps, whereas for small time steps, Krylov-based methods provide better results than Leja-based methods.
\end{abstract}
\section{Introduction}
Exponential integrators are efficient numerical methods for stiff partial differential equations, particularly nonlinear systems (see \cite{HO2010}). These schemes linearize the right-hand side of the problem into linear and nonlinear parts, solving the stiff linear part exactly while treating the nonlinear remainder explicitly. The key challenge is the efficient computation of matrix functions—such as exponentials or $\varphi$-functions of the Jacobian—applied to vectors, which allows for larger time steps compared to standard explicit methods.

For advection-diffusion problems with constant coefficients and simple domains like rectangles with periodic or homogeneous boundary conditions, fast Fourier transform (FFT) and spectral methods are effective for computing these matrix functions (see \cite{Caliari2022, Crouseilles2020, Crouseilles2018, Karle2006}). However, in more general or complex settings, alternative methods are needed. To efficiently approximate the action of the matrix exponential or related functions on vectors, Pad\'{e} approximations or diagonalization are suitable for small systems, whereas for large systems, Krylov subspace methods (see \cite{e6fd4a4b-1a44-3d13-b09d-a560ff78d9f9,10.1145/2168773.2168781}) and Leja interpolation (see \cite{10.1007/s10543-013-0446-0,CALIARI200479}) are commonly used. 
 
Advection-dominated problems are types of partial differential equations where the effect of advection (transport or movement) is much stronger than diffusion (spreading or smoothing). These problems commonly arise in fluid dynamics, transport phenomena, and environmental modeling. In this paper, we consider two common advection-dominated problems include: linear advection–diffusion equations with variable diffusion and two-dimensional compressible isothermal Navier–Stokes systems which describe fluid flow behavior under conditions where advection effects dominate over diffusion.

 Exponential integrators are well-studied for diffusion-dominated problems, where they often outperform explicit and implicit methods (see \cite{doi:10.1137/S1064827595295337,HO2010}). However, they have received less focus for advection-dominated or mixed problems.
We compare the performance of exponential integrators with explicit Runge–Kutta methods by measuring the CPU time required for each simulation, in contrast to constructing a performance model as in \cite{ehoaamm}. For comparison, we evaluate the exponential methods using a 4th-order Krylov approximation in addition to the 2nd-order approximation. We also consider a new regime for numerical testing, in which a localized overpressure generates a propagating shock wave (explosion scenario). For this regime, computations are restricted to the pre-shock stage to ensure the solution remains sufficiently smooth. This approach provides insight into the practical performance, which can be (very) different from theoretical performance models, and the efficiency of each method, complementing the analysis in \cite{ehoaamm}. It also allows us to evaluate how efficiently each method handles different problem regimes. As done in \cite{ehoaamm}, the aim is to examine exponential integrators as general-purpose methods rather than targeting specific problems for which more specialised or efficient approaches might be preferred.
  

The paper is organized as follows. In Section~\ref{sec:methods}, we review the exponential Rosenbrock schemes employed in this study. Section~\ref{testproblem} introduces the test problems used for our experiments. The numerical results and discussions are presented in Sections~\ref{result1d} and \ref{result2d} for one- and two-dimensional problems, respectively. Finally, Section~\ref{conclude} contains the main results of this study.

\section{Exponential Rosenbrock Methods}
\label{sec:methods}
In this study, we are mainly interested in explicit Runge–Kutta methods (see \cite{hairer2008solving,book}) and exponential Rosenbrock methods. We recall the construction of the exponential Rosenbrock methods. Consider the autonomous problem
\begin{equation}\label{2.1}
	u'(t) = F(u(t)), \quad u(t_0) = u_0,
\end{equation}
where \(F\) is a nonlinear function of \(u\). Linearizing \eqref{2.1} gives
\begin{equation}
	u'(t) = J_k u(t) + g_k(u(t)),
\end{equation}
with
\begin{equation}
	J_k = \frac{\partial F}{\partial u}(u_k), \qquad g_k(u(t)) = F(u(t)) - J_k u(t).
\end{equation}

Let \(u_k\) denote the numerical approximation of \(u(t_k)\) at time \(t = t_k\), and let \(\tau\) be the time step size. The exponential Rosenbrock–Euler method \cite{HO2010} is given by
\begin{equation}
	u_{k+1} = u_k + \tau \, \varphi_1(\tau J_k) F(u_k), \qquad 
	\varphi_1(z) = \frac{e^z - 1}{z}.
\end{equation}
This scheme is second-order accurate \cite{HO2010,doi:10.1137/080717717} and requires computing the action of the matrix function \(\varphi_1\) on a vector.  We also take into account higher-order methods (see \cite{doi:10.1137/080717717}).  
Recently, efficient two-stage fourth-order methods for time-dependent PDEs have been proposed \cite{LUAN201791}. In this paper, we employ the \(\mathtt{exprb42}\) scheme:
\begin{equation}
\begin{aligned}
	U_{k2} &= u_k + \frac{3}{4} \, \tau \, \varphi_1\Big(\frac{3}{4} \tau J_k\Big) F(u_k), \\[1ex]
	u_{k+1} &= u_k + \tau \, \varphi_1(\tau J_k) F(u_k) + \frac{32}{9} \, \tau \, \varphi_3(\tau J_k) \Big(g(U_{k2}) - g(u_k)\Big),
\end{aligned}
\end{equation}
where
\begin{equation}
	\varphi_3(z) = \frac{e^z - 1 - z - z^2/2}{z^3}.
\end{equation}
This method is computationally attractive because it requires only one internal stage per time step and involves three evaluations of matrix functions, which is fewer than the related \(\mathtt{exprb43}\) scheme \cite{HO2010}.

\section{Model problems}\label{testproblem}
First, we study the linear advection–diffusion equation, where diffusion is relatively weak compared to advection, in Section~\ref{subsec: linear}. Next, the nonlinear compressible isothermal Navier-Stokes problem, which shows strong advection-dominated behavior, is introduced in Section~\ref{discussNavierstokes}.

\subsection{Linear advection-diffusion problem}
\label{subsec: linear}

We consider the linear advection–diffusion problem 
\begin{equation}\label{problem1}
	\begin{aligned}
		& \partial_t  u(t,x) - \kappa(x) \partial_{xx} u (t,x) + \partial_x u(t,x)= 0,\quad (t,x)\in(0,1)\times[0,1], \\
		& u(t,0)= u(t,1)=0,\\
		& u(0,x) = x(1-x),
	\end{aligned}	
\end{equation}
with homogeneous Dirichlet boundary conditions. Here \( \kappa\) is the diffusion coefficient. For spatial discretization of problem \eqref{problem1}, we use a uniform grid with $n = 1599$ interior grid points with a mesh size $h = 1/(n+1) = 1/1600$. We approximate both the diffusion operator $-\kappa(x)\partial_{xx}$ and the advection operator $\partial_x$ using standard second-order finite differences, resulting in discrete matrices $A_h$ and $B_h$, respectively.

\begin{figure}[h!]
\subfigure{\includegraphics[width=0.51\textwidth]{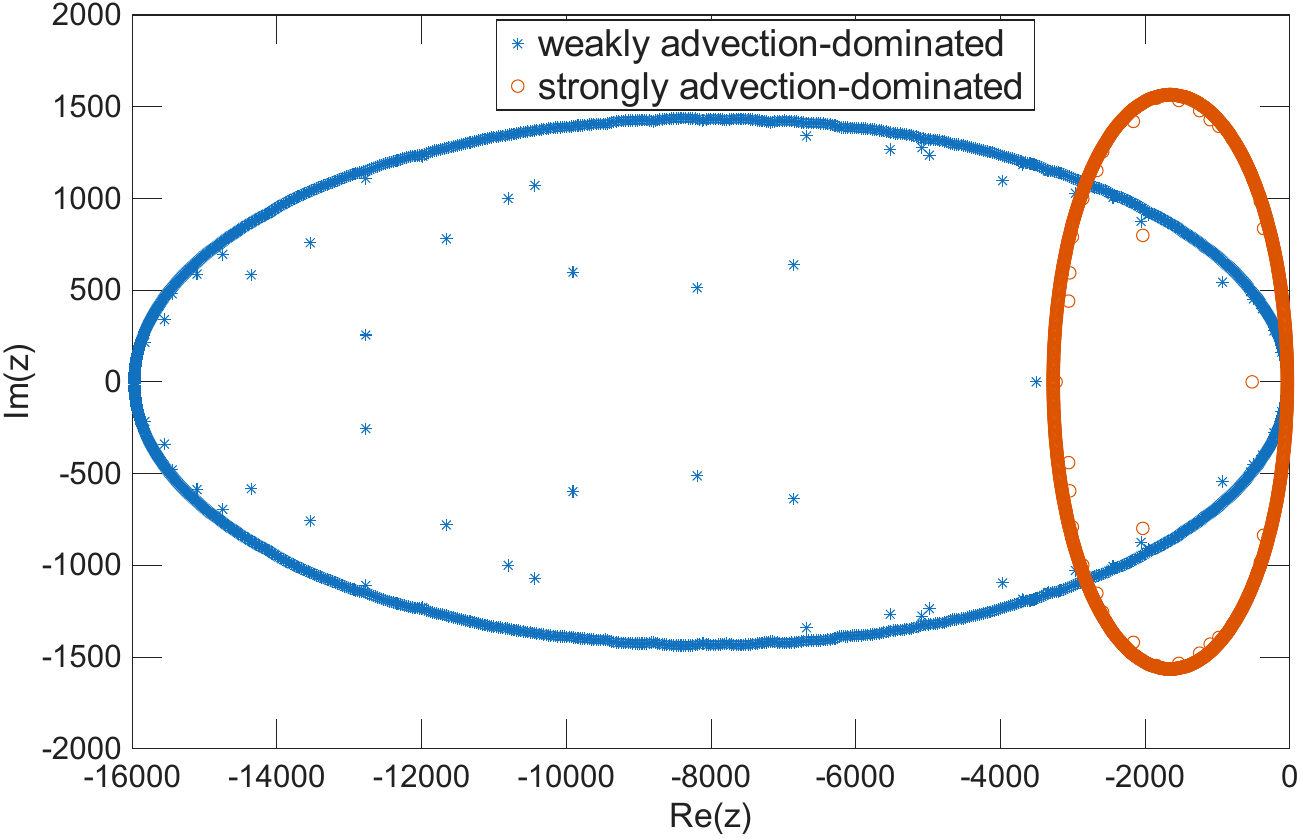}}
	\subfigure{\includegraphics[width=0.51\textwidth]{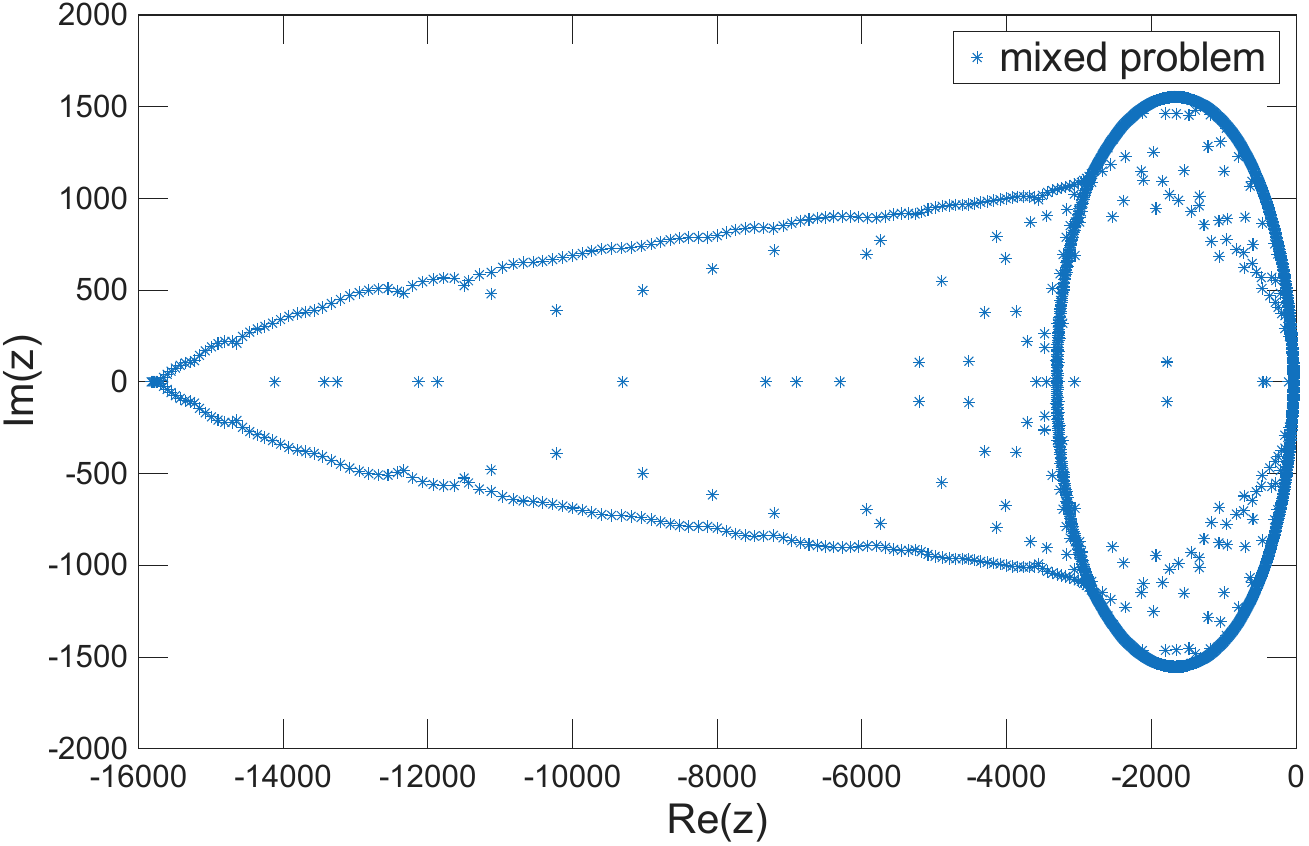}}\label{fig4c}			
	\caption{Numerical spectrum of the matrix $-(A_h + B_h)$ for problem \eqref{problem1} for different values of $\kappa$. The left panel corresponds to $\kappa = 1/640$,  representing a weakly advection-dominated case, while $\kappa = 1/3100$ corresponds to a strongly advection-dominated case. The figure on the right shows the spectrum for a spatially varying $\kappa$, illustrating regions that are weakly or strongly advection-dominated depending on the location.} 
	\label{spectrum1}
\end{figure}

All computations are performed in the advection-dominated regime, characterized by small values of $\kappa \ll 1$. We consider two cases:
\begin{enumerate}
	\item \emph{weakly advection-dominated} regime, where the stability restriction (CFL condition) imposed by diffusion is more restrictive than that imposed by advection.
	\item \emph{strongly advection-dominated} regime, where the CFL restrictions for diffusion and advection are approximately equal.
\end{enumerate}
Figure~\ref{spectrum1} displays the numerical spectrum (by using the MATLAB routine $\mathtt{eigs}$) of the matrix $-(A_h+B_h)$.

\subsection{2D Compressible isothermal Navier--Stokes problem}\label{discussNavierstokes}
The isothermal Navier--Stokes system in two spatial dimensions is given by
\begin{equation}\label{eq41}
	\begin{aligned}
		&\partial_t \rho + \partial_x (\rho u) + \partial_y (\rho v)  = 0,	\\
		&\partial_{t} u + u \partial_{x} u+v \partial_{y} u + \tfrac{1}{\rho}\partial_{x} \rho =  \nu\left(\partial_{x x} u+\partial_{y y} u\right),	\\
		&\partial_{t} v+u \partial_{x} v+v \partial_{y} v +\tfrac{1}{\rho} \partial_{y} \rho =\nu\left(\partial_{x x} v+\partial_{y y} v\right).
	\end{aligned}
\end{equation}
The first equation represents the conservation of mass, while the latter two correspond to the momentum equations. 
The parameter $\nu$ is a positive constant. 
We consider problem \eqref{eq41} on the domain $\Omega = [0,1]^2$ with periodic boundary conditions over the time interval $t \in [0,12]$. 
The initial values are given by the explosion and shear flow initial conditions.

\begin{figure}[h!]
	\subfigure{\includegraphics[width=0.5\textwidth]{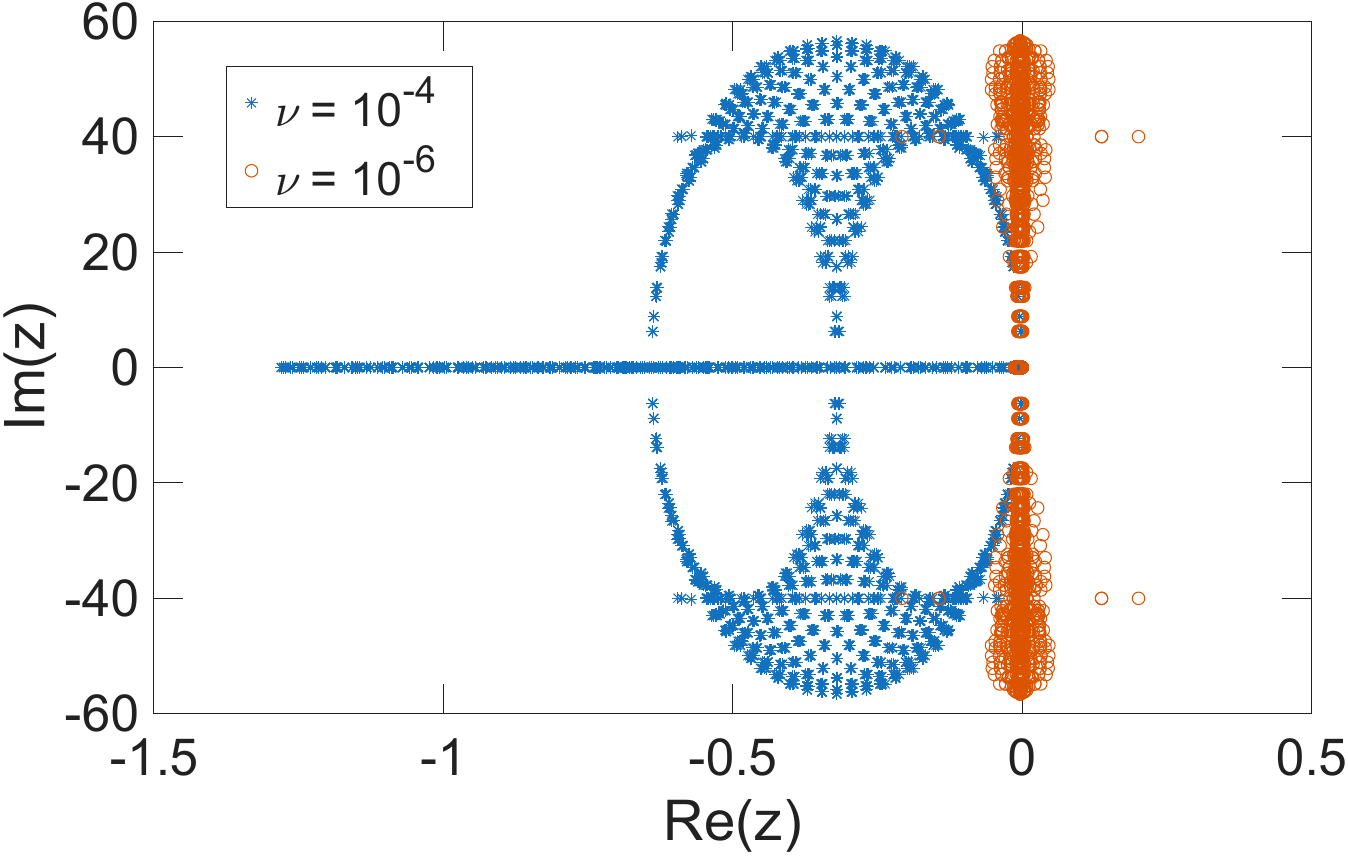}}
	\subfigure{\includegraphics[width=0.5\textwidth]{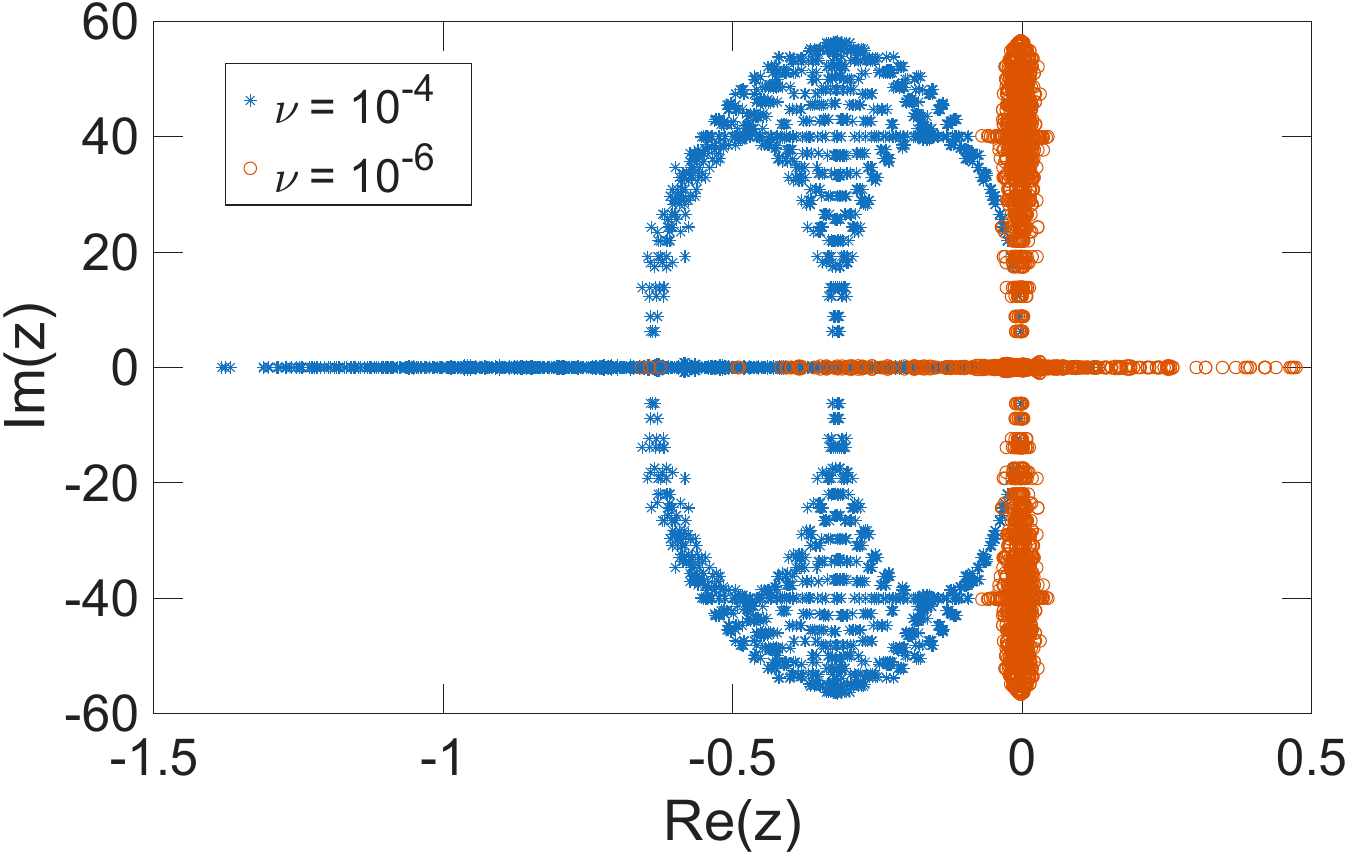}}	
	\caption{The numerical spectra of the Jacobians evaluated at the initial value (left figure) and at the final value (right figure) of a numerical solution of the two-dimensional compressible isothermal Navier--Stokes problem \eqref{eq41} with diffusion coefficients of \(\nu = 10^{-4}\) and \(\nu = 10^{-6}\), respectively, indicating strongly advection-dominated scenarios.}
	\label{spectrum3}
\end{figure}

We reformulate problem \eqref{eq41} as follows
\begin{equation}
	\label{1729}
	U^{\prime}(t)=F(U(t)), \quad U(0)=U_{0} ,
\end{equation}
where
\begin{equation}\label{eq42}
	U = \left[\begin{array}{l}
		\rho \\
		u \\
		v
	\end{array}\right], \quad F(U) = \left[\begin{array}{c}
		F_1(U) \\
		F_2(U) \\
		F_3(U)
	\end{array}\right] =
	\left[\begin{array}{c}
		-\partial_x (\rho u) - \partial_y (\rho v) \\
		- u \partial_{x} u - v \partial_{y} u - \frac{1}{\rho} \partial_{x} \rho + \nu (\partial_{xx}u + \partial_{yy}u) \\
		-u \partial_{x} v-v \partial_{y} v - \frac{1}{\rho} \partial_{y} \rho +\nu\left(\partial_{x x} v+\partial_{y y} v\right)
	\end{array}\right].
\end{equation}
The parameter $\nu$ is crutial in numerical experiments involving diffusion processes, as it directly controls the strength of diffusion in the system. Varying $\nu$ allows us to investigate how diffusion influences solution behavior, stability, and the transition between different dynamical regimes. In \cite{ehoaamm} provides detailed procedures on constructing the Jacobian matrix and the right-hand side vector $F(U)$ in convection-diffusion problems. We provide the numerical spectrum of the Jacobian matrix in the explosion case (see \cite{ehoaamm} for the shear-flow case), using 40 degrees of freedom, as illustrated in Figure~\ref{spectrum3}.

\section{Numerical results for linear advection-dominated problems}\label{result1d}
This section focuses on evaluating the efficiency and accuracy of various numerical methods for advection-dominated problems, using a fixed $n=1599$ interior grid points.  The exact solution is used for the linear problem. The $L^2$ error of each method at the final time \(T\) is then calculated and reported along with the corresponding CPU time.

Henceforth, all tests are performed in MATLAB R2025a on a standard laptop equipped
with an Intel Core i7-13620H processor with 10 physical cores and 16 GB of RAM. Parallelization and GPU acceleration are not used in this study.

\subsection{Weakly advection-dominated problem}
Setting $\kappa = \frac{1}{640}$ places problem \eqref{problem1} in a weakly advection-dominated regime. In this context, both Krylov subspace and Leja interpolation methods demonstrate significant advantages over  explicit Runge--Kutta schemes. Both Krylov and Leja methods require us to specify a time step size and a target accuracy (tolerance). In our numerical experiments, the tolerances $10^{-4}$ and $10^{-7}$ are used. We choose a maximum time step size of $\tau = \frac{1}{48}$ for Leja and $\tau = \frac{1}{12}$ for Krylov. The numerical results are shown in Figure~\ref{fig1}. It is observed that 
 Krylov and Leja methods permit much larger time steps compared to explicit methods. In particular, Krylov methods can allow time steps up to 750 times larger, and Leja methods up to 188 times larger, than those permitted by RK2 and RK4 in advection-dominated regimes highlighting the substantial efficiency gains possible with these exponential integrators. 

Both Krylov and Leja methods are more computationally efficient than explicit Runge--Kutta schemes, but their relative efficiency depends on the time step size and underlying operations. The Leja method relies solely on matrix-vector multiplications, making it highly efficient for large time steps and large-scale problems. Krylov methods use the Arnoldi process to build an orthonormal basis, which becomes increasingly expensive for large time step. Large time steps require a longer recurrence (more basis vectors), leading to higher computational and memory costs. For smaller time steps, the required Krylov subspace is smaller, reducing the cost of orthogonalization and making Krylov methods more competitive or even superior in some integrator-specific scenarios.

\begin{figure}[h!]
	\subfigure[ $\tau = \frac{1}{48}$ ]{\includegraphics[width=0.49\textwidth]{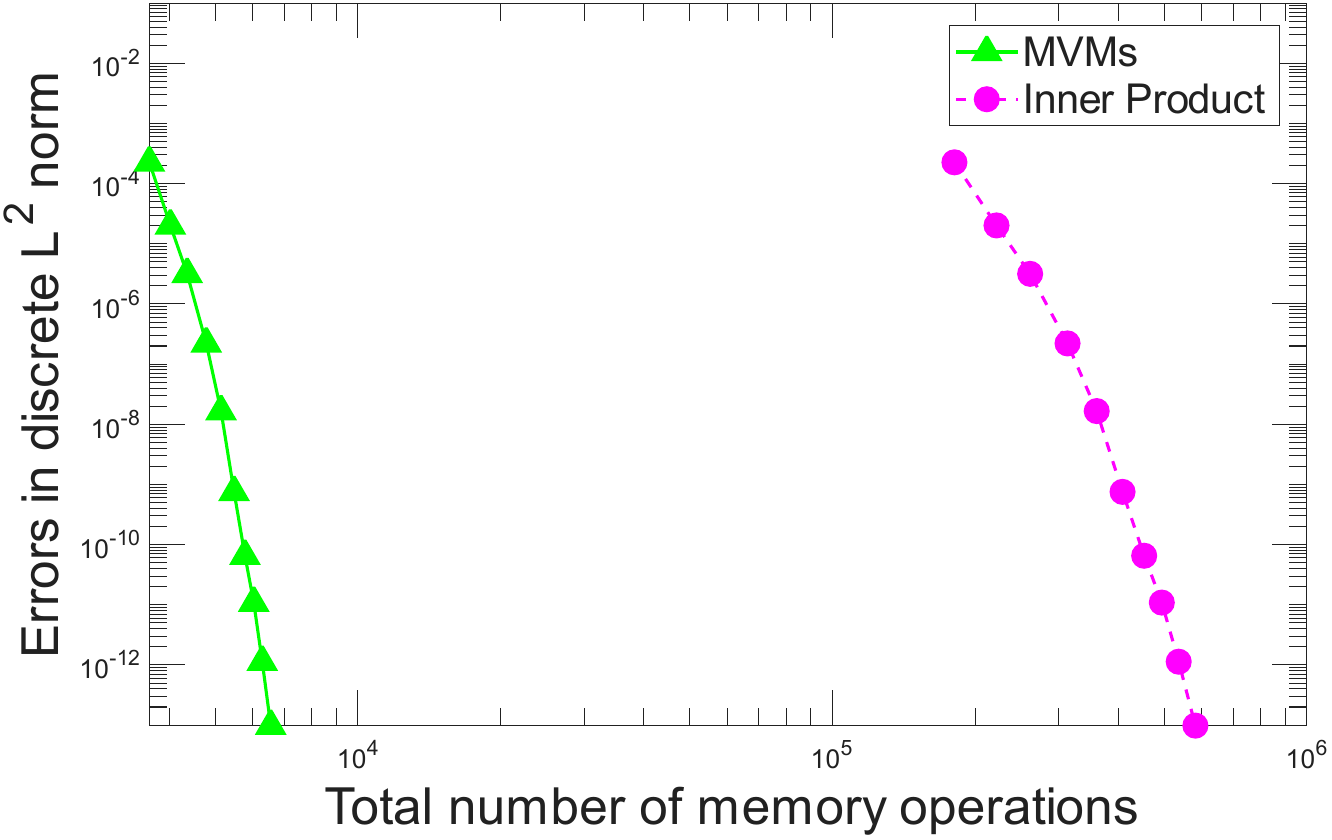}\label{fig1e}}
	\subfigure[ $\tau = \frac{1}{48}$]{\includegraphics[width=0.49\textwidth]{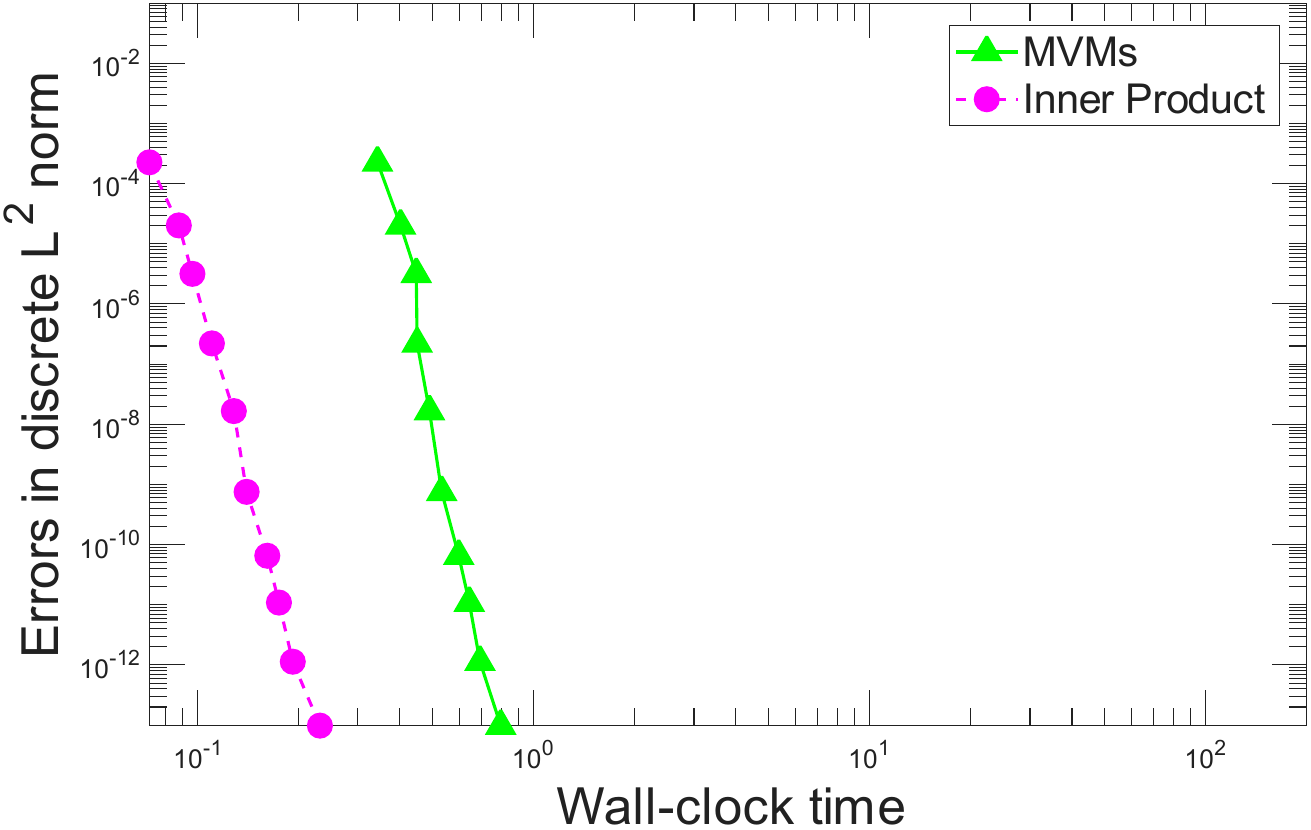}\label{fig1f}}		
	\caption{The numerical results for the weakly advection-dominated problem \eqref{problem1} with $\kappa = \frac{1}{640}$. Figure~(a) shows the total memory operations of MVMs as a function of the $L^2$ error, while Figure~(b) depicts the achieved accuracy of the considered methods as a function of wall-clock time.}
\end{figure}

\begin{figure}[h!]
		\subfigure[]{\includegraphics[width=0.51\textwidth]{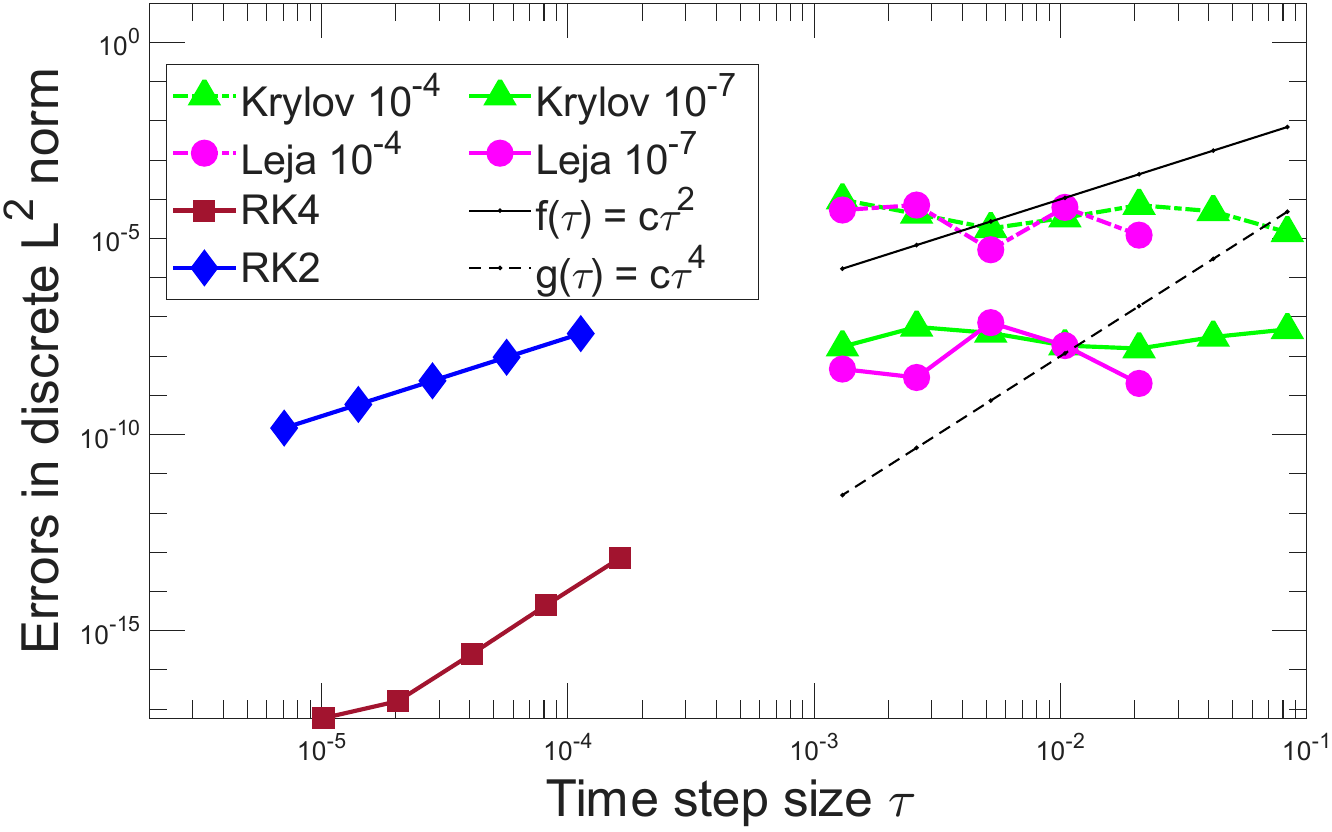}\label{fig1a}}
		\subfigure[]{\includegraphics[width=0.51\textwidth]{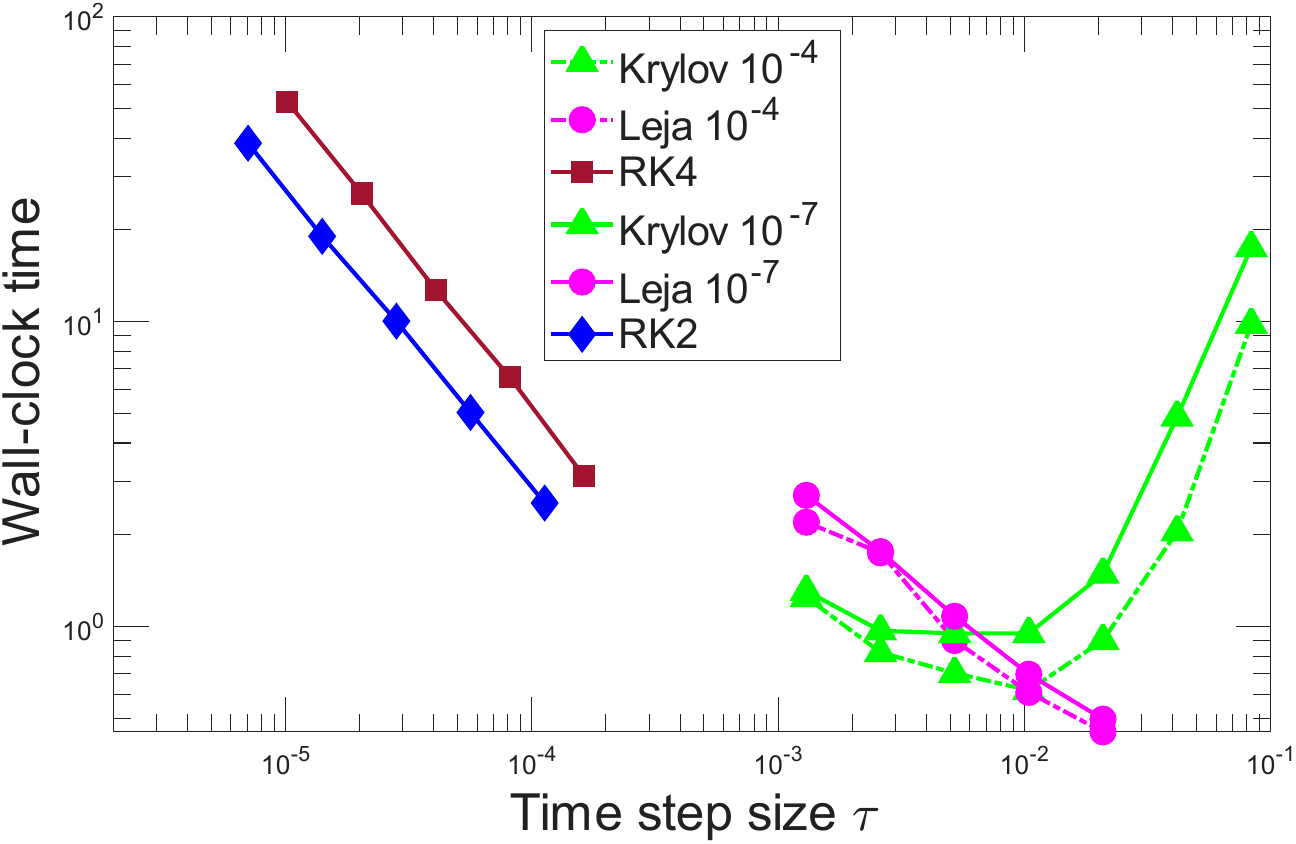}\label{fig1b}}	
		\subfigure[ $\tau = \frac{1}{48}$ ]{\includegraphics[width=0.51\textwidth]{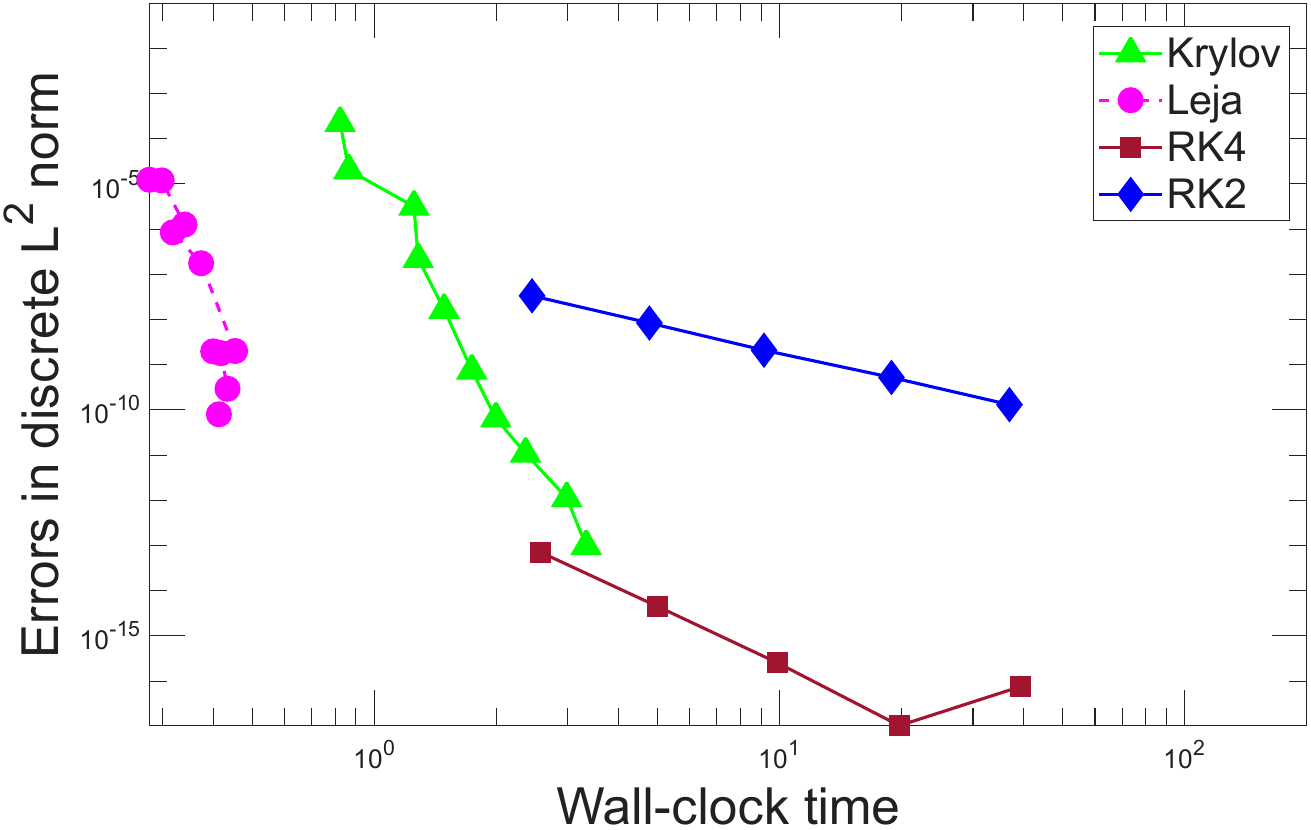}\label{fig1c}}
		\subfigure[ $\tau = \frac{1}{768}$]{\includegraphics[width=0.51\textwidth]{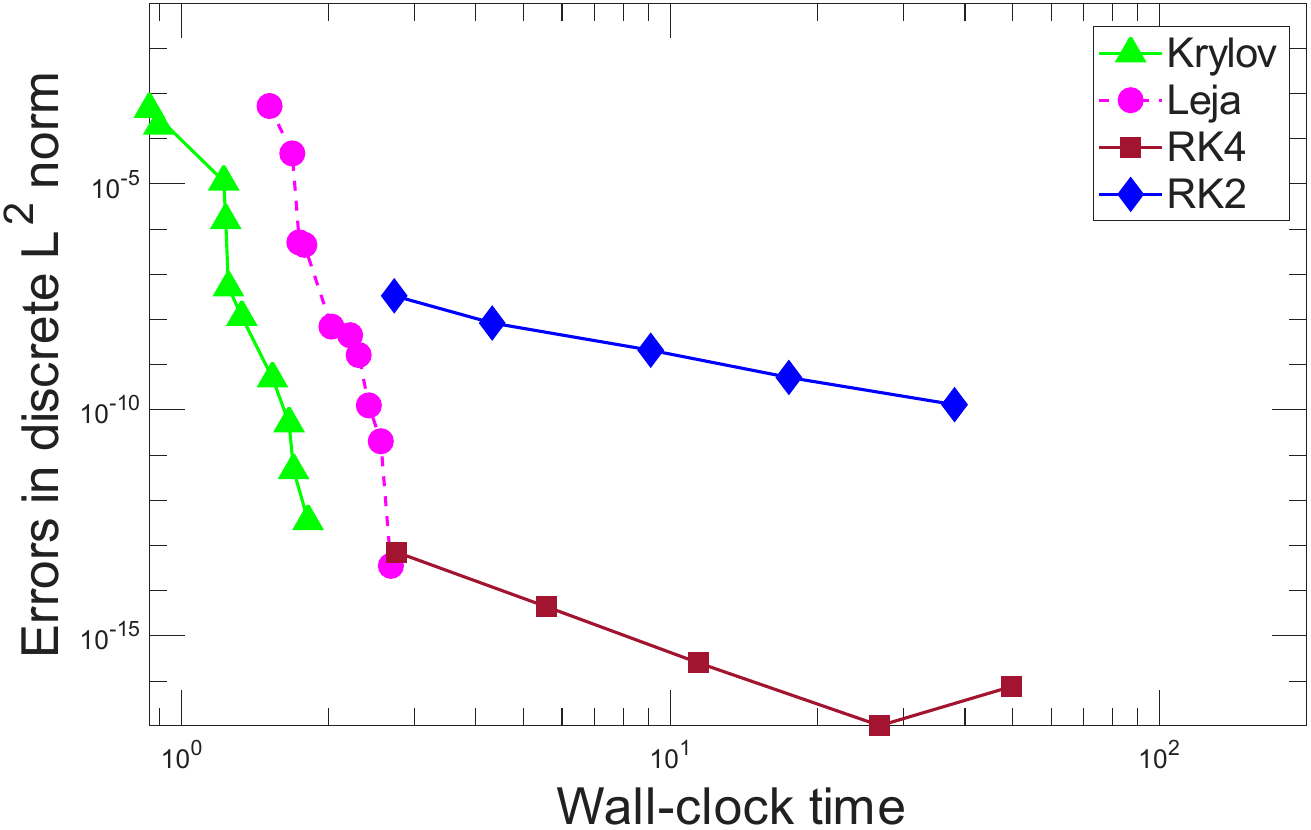}\label{fig1d}}		
		\caption{The numerical results for the weakly advection-dominated problem \eqref{problem1} with $\kappa = \frac{1}{640}$ are shown in the top two figures of Figure~\ref{fig1}. These figures display the accuracy achieved and the computational cost for different methods evaluating $\exp(- (A_h + B_h)) u_0$ as a function of the time step size. The tolerances are set so that the exponential integrators reach accuracies of $10^{-4}$ (dashed-dotted line) and $10^{-7}$ (solid line), respectively. The bottom two panels show the accuracy of the methods as a function of computational cost for two selected time step sizes $\tau$.
			 \label{fig1}}
	\end{figure}

Moreover, Figure~\ref{fig1} shows that Krylov methods can even outperform the Leja method in small step sizes. This is not in contradiction with the results in \cite{ehoaamm}, which indicate that Krylov methods generally require more memory operations than Leja. Although the Krylov subspace method involves additional inner products due to orthogonalization, these operations are relatively inexpensive compared to the operator applications required by both methods. Each inner product only requires sequential access to two vectors, benefiting from good cache locality. In contrast, the Leja method performs matrix-vector actions, which dominate the computational cost. These evaluations typically involve more memory traffic and arithmetic per step than the inner products in Krylov. As a result, Leja interpolation is often more memory-bound, whereas Krylov methods can achieve lower wall-clock times in practice due to more efficient data reuse within the subspace basis. Although there are far fewer matrix–vector multiplications (MVMs) than inner products, the execution time of inner products is still significantly faster than that of MVMs. 
These observations are confirmed in Figure \ref{fig1e}, Figure \ref{fig1f}, Figure \ref{fig1c}, and Figure \ref{fig1d}. 

At a moderate accuracy level of $10^{-4}$, both Leja and Krylov methods efficiently approximate the action of the exponential function, while explicit Runge--Kutta schemes fail to achieve this. Even at a stricter tolerance of $10^{-7}$, Krylov and Leja methods still outperform classical integrators like RK2 and RK4. Figure~\ref{fig1} shows that exponential integrators can reach much higher accuracy without significantly increasing computational cost, unlike explicit methods. Notably, the Krylov and Leja method can achieve any desired accuracy using fewer computational resources compared to RK2 and RK4.

\subsection{Strongly advection-dominated problem}
Setting $\kappa = \frac{1}{3100}$ yields a grid Péclet number of approximately $ \frac{3100}{1600} \approx 1.93$, placing problem \eqref{problem1} in a strongly advection-dominated regime. This value is close to the threshold where centered finite difference discretizations remain non-oscillatory for the given spatial resolution (see \cite{hundsdorfer2013numerical}), as high P\'eclet  numbers are known to induce spurious oscillations in standard schemes. For both the Krylov and Leja time integration methods, we choose a maximum time step of $\frac{1}{102}$.

\begin{figure}[t!]
		\subfigure[]{\includegraphics[width=0.51\textwidth]{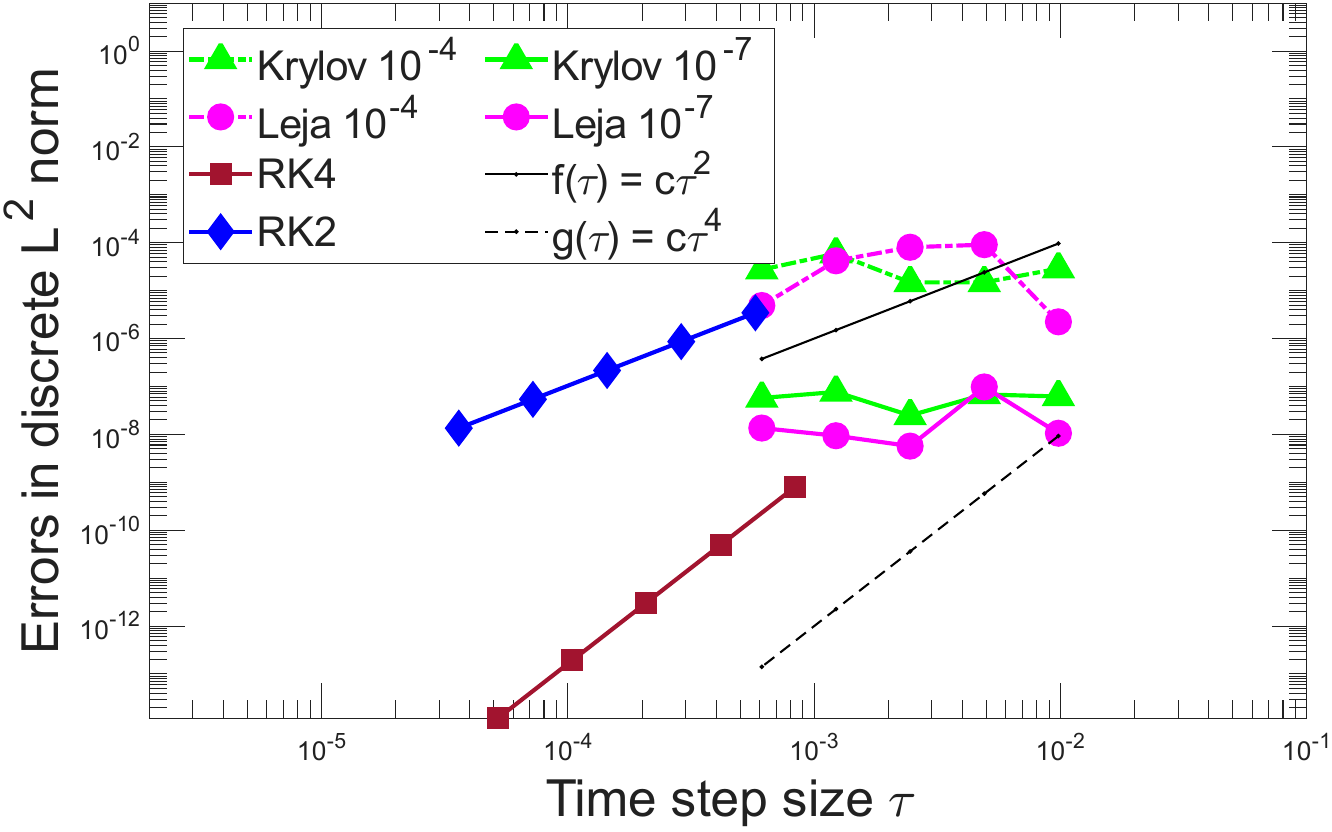} \label{fig2a}}
		\subfigure[]{\includegraphics[width=0.51\textwidth]{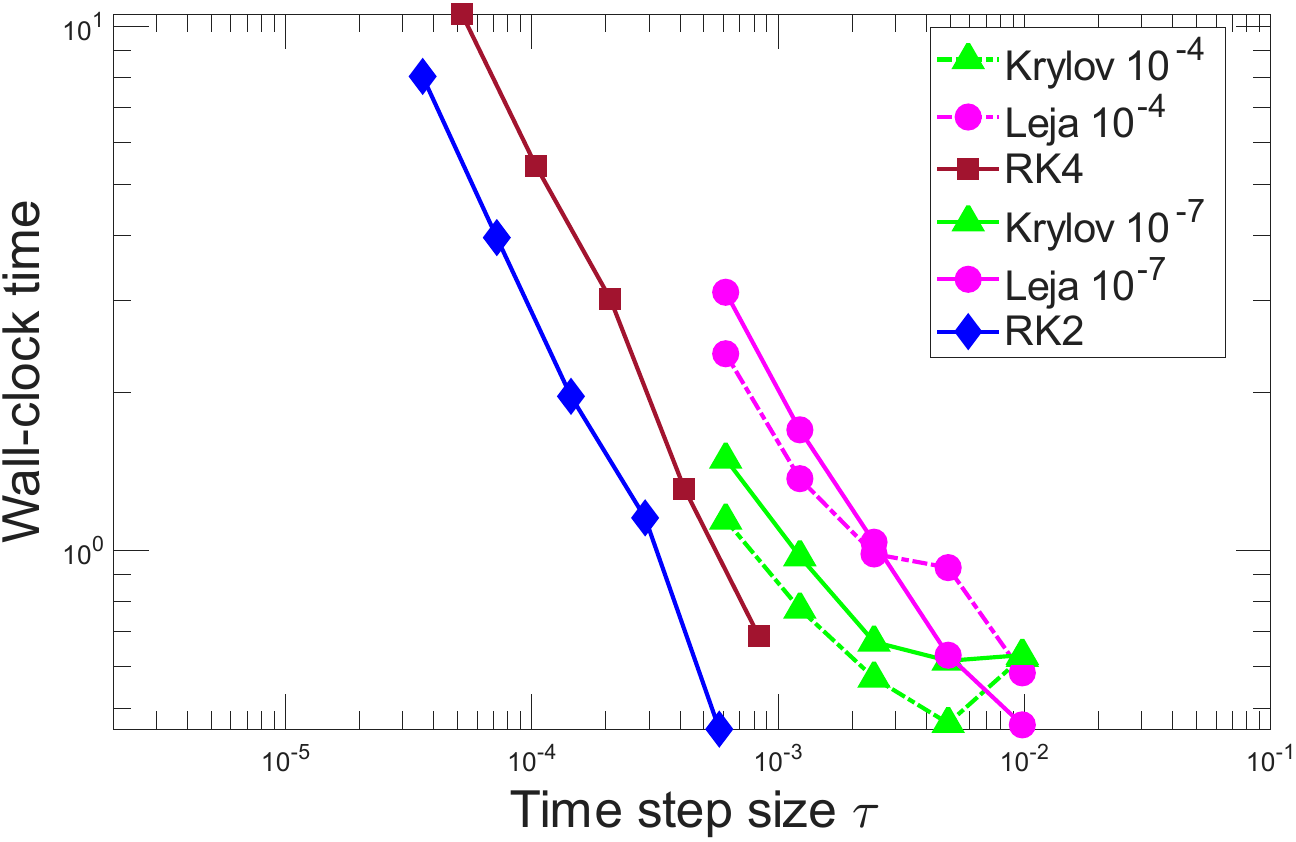} \label{fig2b}}	
		\subfigure[$\tau = \frac{1}{102}$ ]{\includegraphics[width=0.51\textwidth]{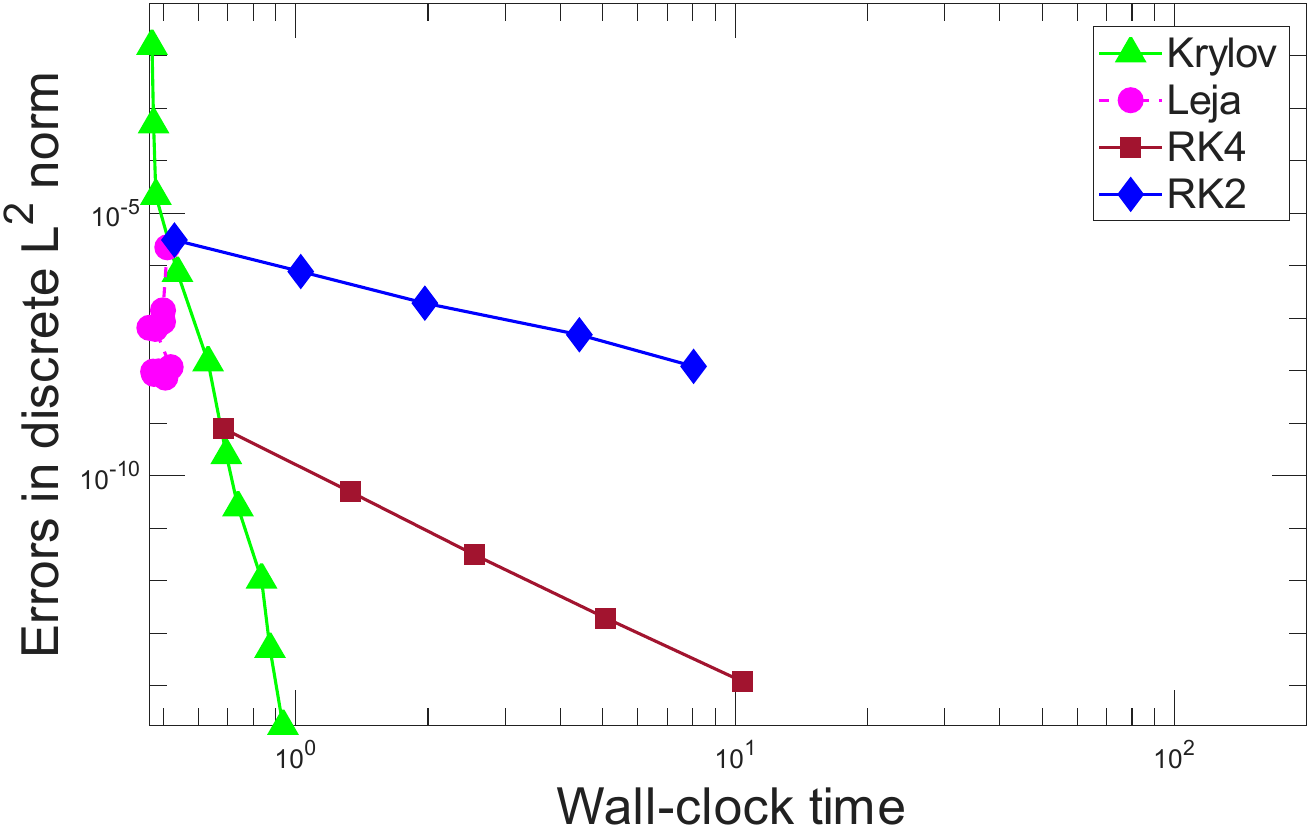} \label{fig2c}}
		\subfigure[ $\tau = \frac{1}{1632}$ ]{\includegraphics[width=0.51\textwidth]{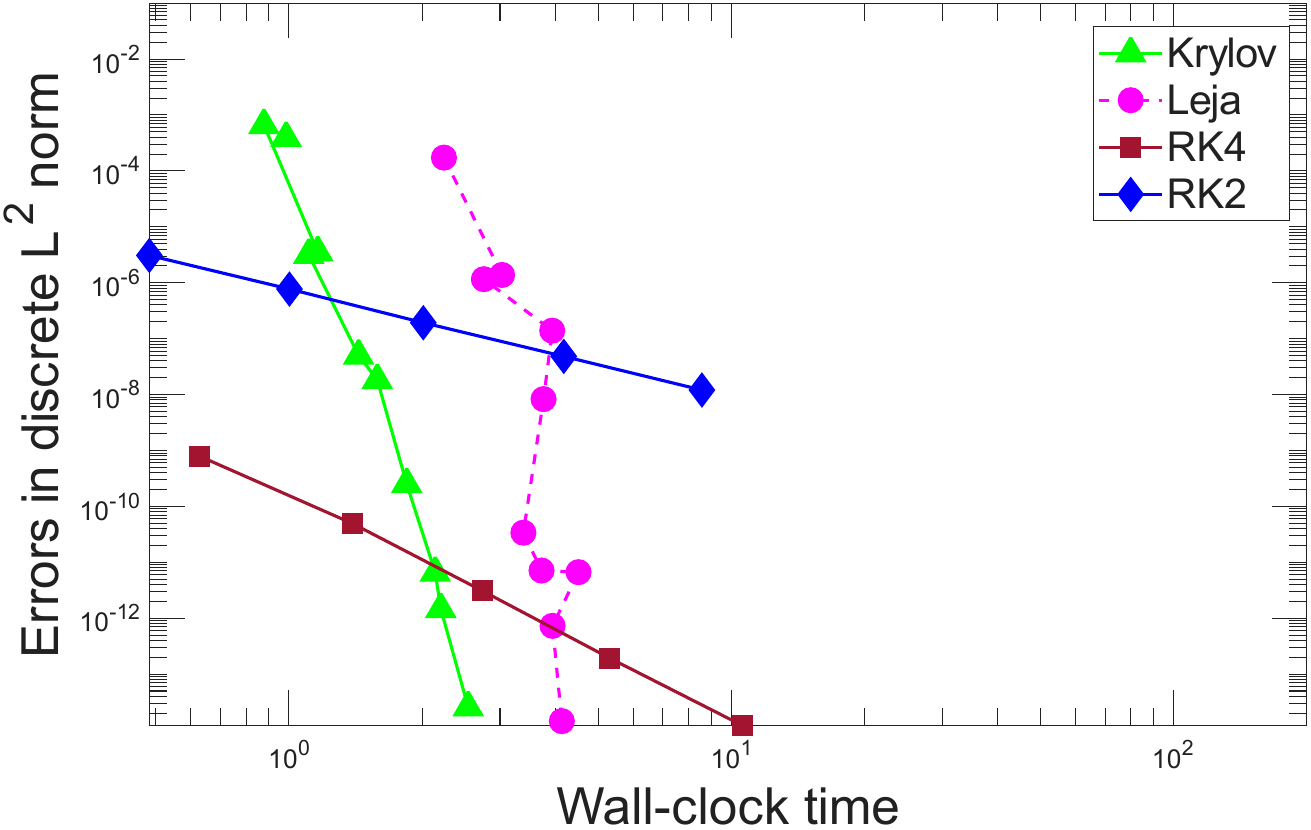} \label{fig2d}}					
		\caption{The top two panels present the numerical results for the strongly advection-dominated problem \eqref{problem1} with $\kappa = \frac{1}{3100}$. These figures show the achieved accuracy and computational cost of the considered methods for evaluating $\exp(- (A_h + B_h)) u_0$ as a function of the time step size. The tolerances are set so that the exponential integrators reach target accuracies of
 $10^{-4}$ (dashed-dotted line) and $10^{-7}$ (solid line), respectively. The bottom two panels illustrate how the accuracy of these methods varies with computational cost for two selected time step values $\tau$.
			 \label{fig2}}
	\end{figure}

Figure~\ref{fig2} shows the results of our numerical experiments. It is observed that exponential integrators allow much larger time steps than RK2 and RK4 schemes- about 17 times larger than RK2 and 12 times larger than RK4. Both the Krylov and Leja methods demonstrate superior computational efficiency compared to explicit Runge--Kutta methods. The Leja method becomes more expensive when a small time step is used, as observed in the weakly advection-dominated case, whereas the Krylov method can benefit from a small time step size.  Overall, the exponential integrators achieve substantially higher accuracy without a notable increase in computational cost, in contrast with the explicit Runge--Kutta schemes.

\subsection{Mixed problem}
\label{mixproblem}
In this section, we study scenarios where advection dominates differently across the spatial domain. For this purpose, we consider problem \eqref{problem1} with the diffusion coefficient \(\kappa\) is set to 
\begin{equation}\label{eq12}
	\kappa(x) = \frac{1}{1067.2} + \frac{1}{1612.5} \tanh(20x - 16),
\end{equation}
so that $\kappa(x)$ varies smoothly from a strongly advection-dominated case at $x = 0$ to a weakly advection-dominated regime at $x = 1$(see Figure~\ref{fig4cc}). In our experiments, the maximum time step sizes for matrix exponential evaluation are $\tau = \frac{1}{48}$ with the Krylov method and $\tau = \frac{1}{96}$ with the Leja method.

The numerical results in Figure~\ref{fig4} demonstrate that both Krylov and Leja exponential integrators chieve high accuracy with significantly lower computational cost than explicit Runge--Kutta methods. This superiority is particularly evident under stringent tolerances, where explicit methods require much smaller time steps and increased computational effort. In this mixed regime problem, the Krylov and Leja approaches exhibit comparable performance.

	\begin{figure}[t!]
		\begin{center}
			\subfigure {\includegraphics[width=0.51\textwidth]{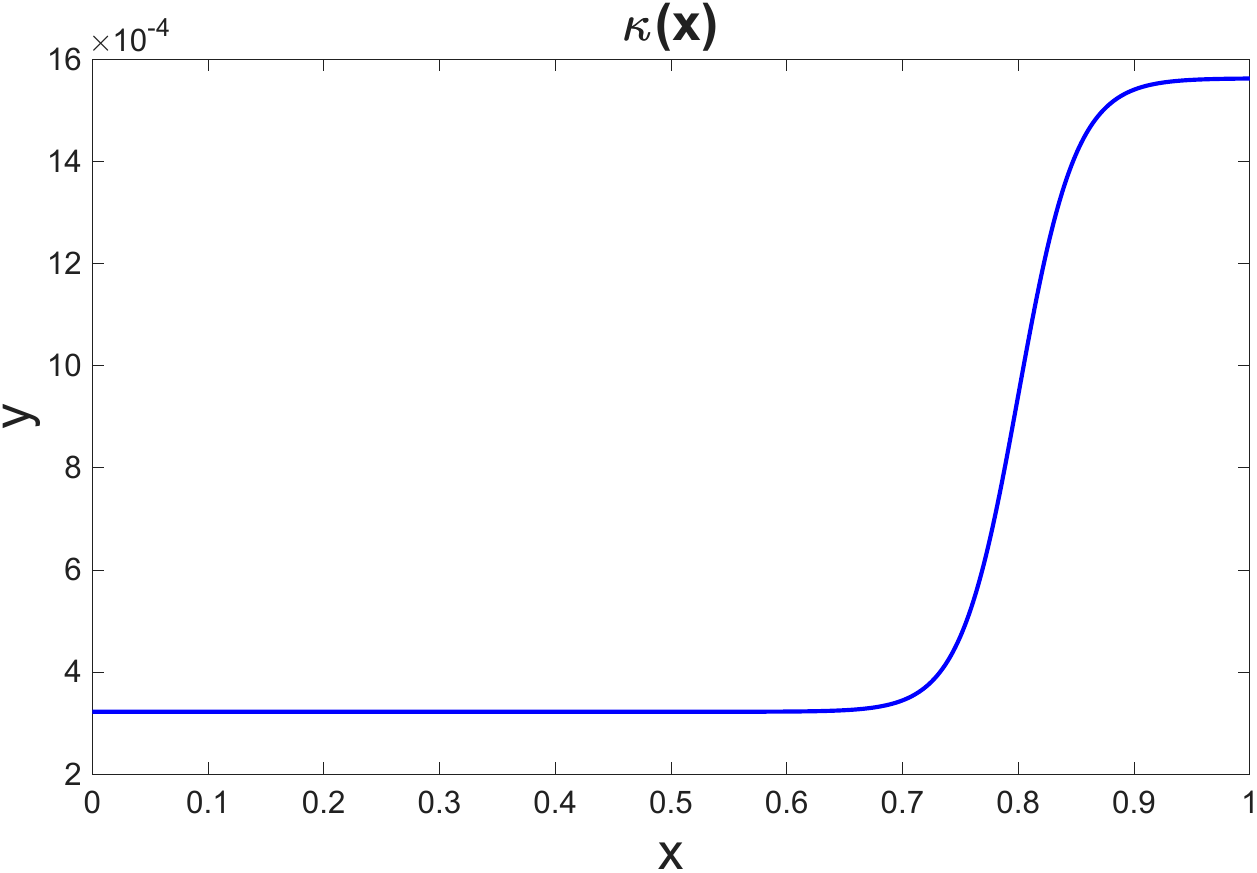}}
		\end{center}							
		\caption{Illustration of the function $\kappa$ defined in~\eqref{eq12}. }\label{fig4e}
	\end{figure}
	
	\begin{figure}[t!]
		\subfigure[]{\includegraphics[width=0.51\textwidth]{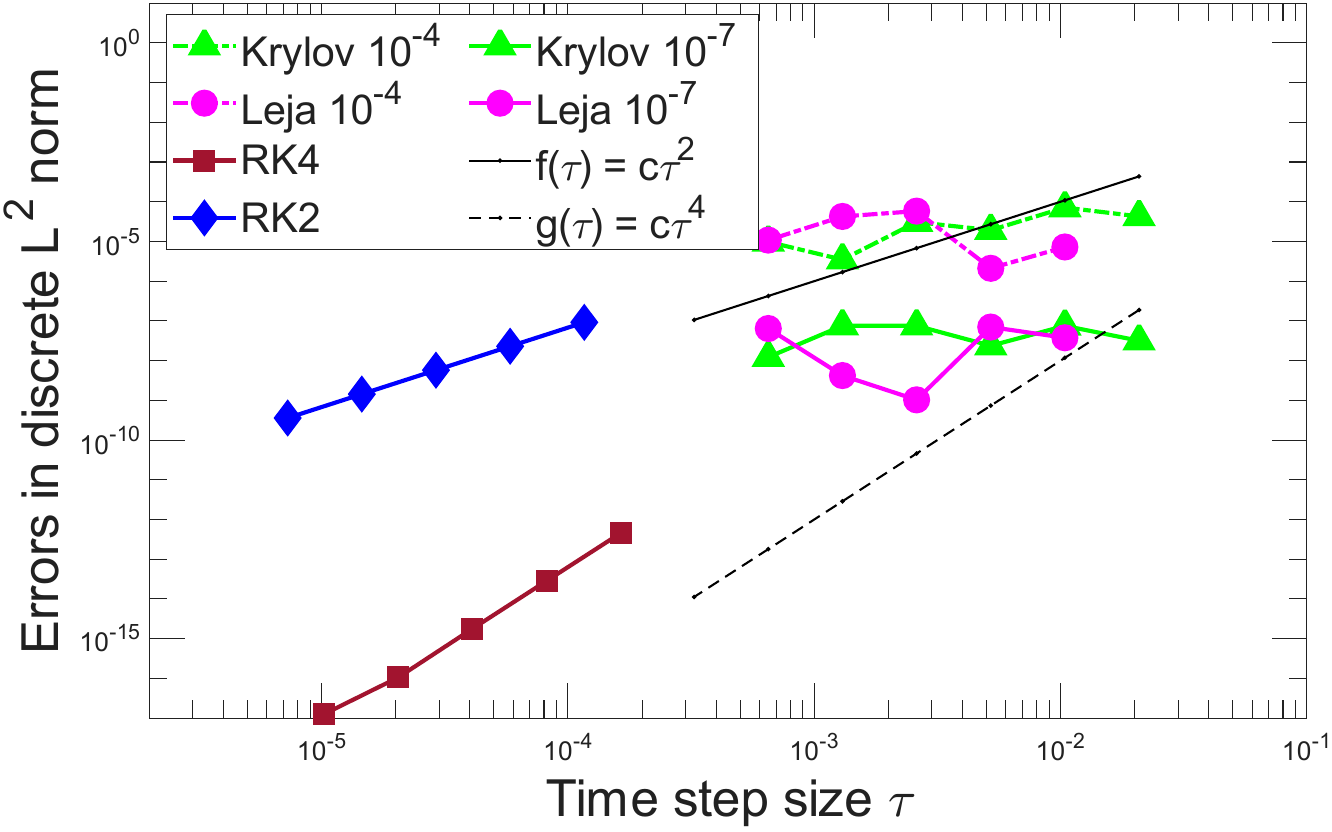} \label{fig4a}}
		\subfigure[]{\includegraphics[width=0.51\textwidth]{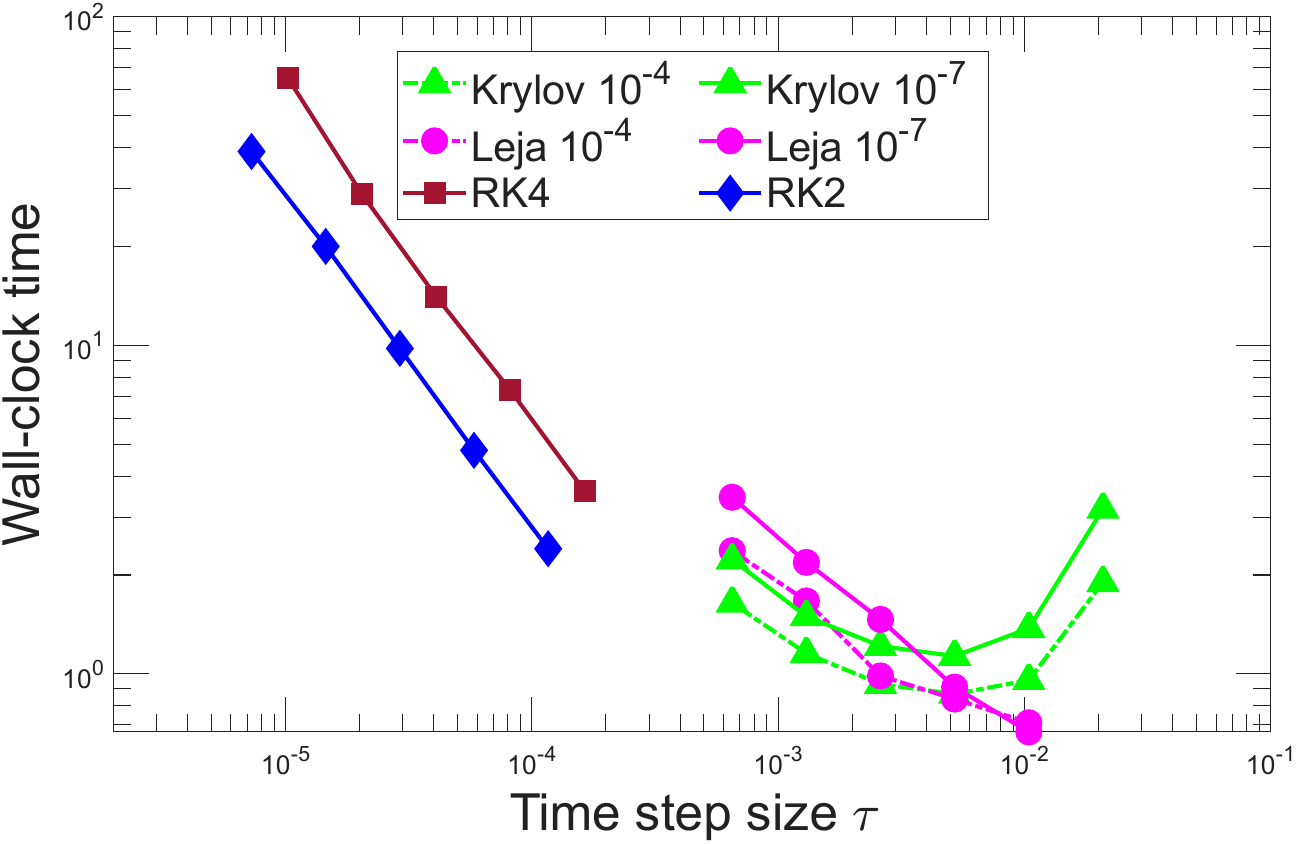} \label{fig4b}}	
		\subfigure[$\tau = \frac{1}{96}$ ] {\includegraphics[width=0.51\textwidth]{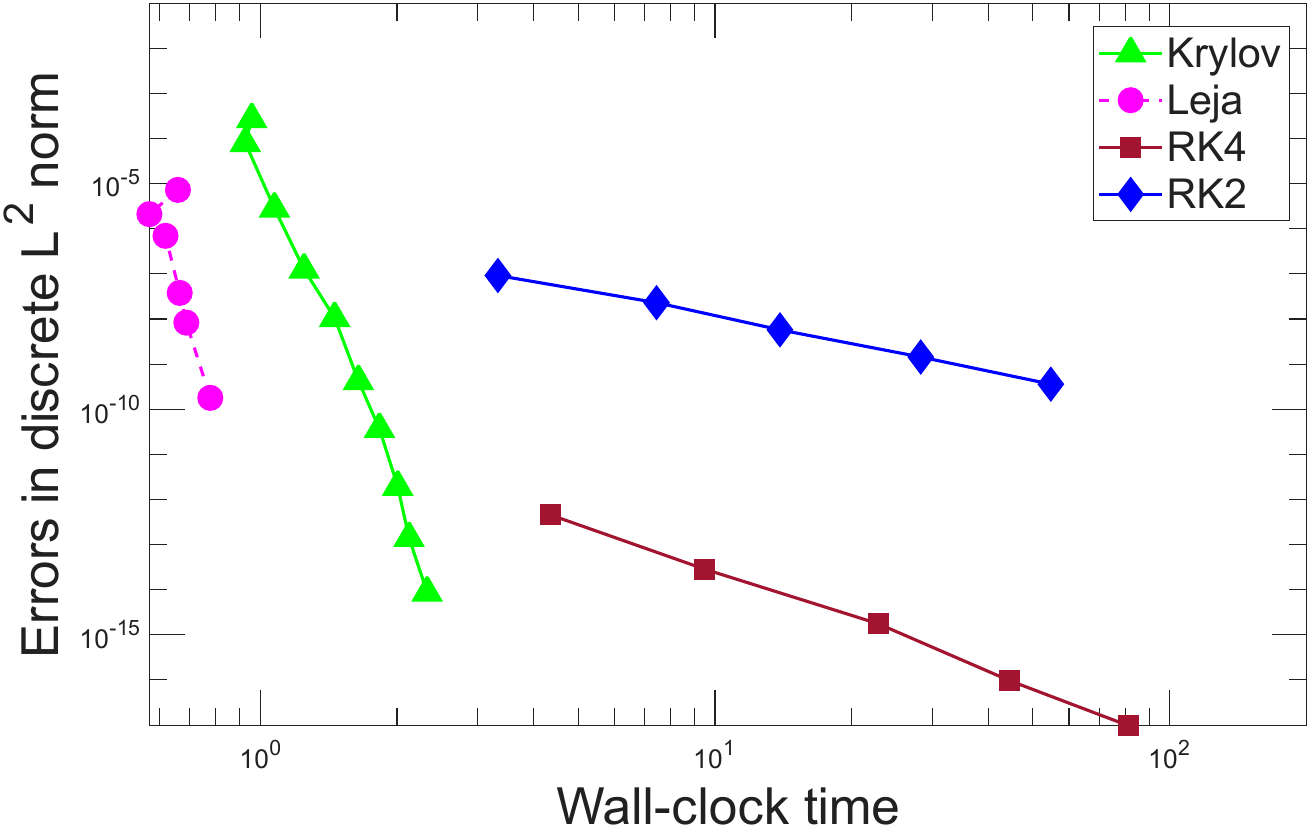} \label{fig4cc}}
		\subfigure[$\tau = \frac{1}{384}$ ]{\includegraphics[width=0.51\textwidth]{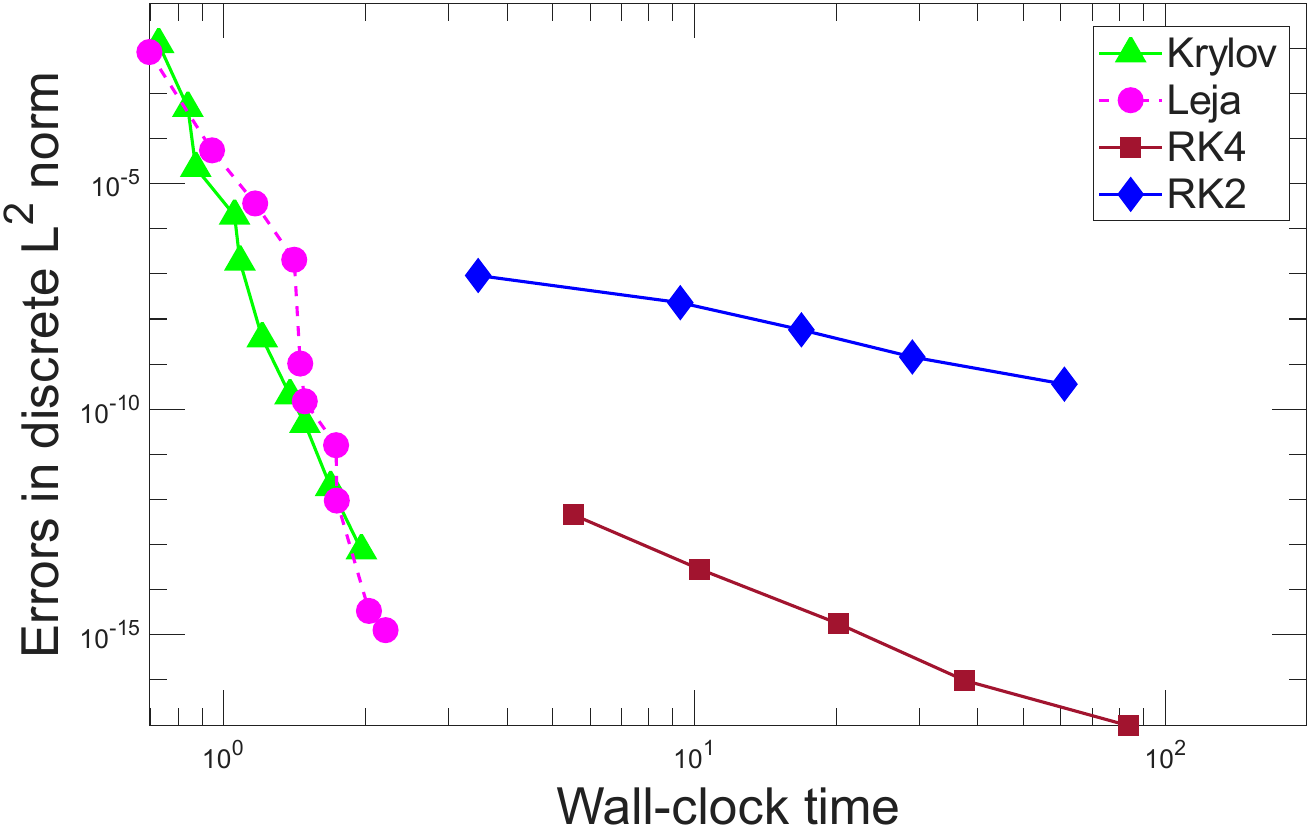}\label{fig4d}}
		\caption{
		The numerical results for the mixed problem \eqref{problem1} with \(\kappa\) in \eqref{eq12}, appear in the upper two panels. These plots display the accuracy and computational cost of the tested methods for evaluating $\exp(- (A_h + B_h)) u_0$ as a function of the time step size. The tolerances are chosen such that the exponential integrators attain accuracies of $10^{-4}$ (dashed-dotted line) and $10^{-7}$ (solid line), respectively. The lower two panels depict the achieved accuracy of the methods as a function of computational cost for two selected values of $\tau$.
			}\label{fig4}
	\end{figure}

\section{Numerical results for the compressible isothermal Navier--Stokes problem}	
\label{result2d}
Encouraged by the results obtained for the linear problem, we now turn to a more challenging nonlinear test case to further evaluate the efficiency of exponential integrators. We study exponential integrators for solving the two-dimensional compressible isothermal Navier–Stokes equations. This test examines whether the advantages observed in linear cases extend to nonlinear fluid dynamics problems.

\subsection{Explosion} 
\begin{figure}[t!]
	\subfigure[$\rho(0.4,x,y)$]{\includegraphics[width=0.51\textwidth]{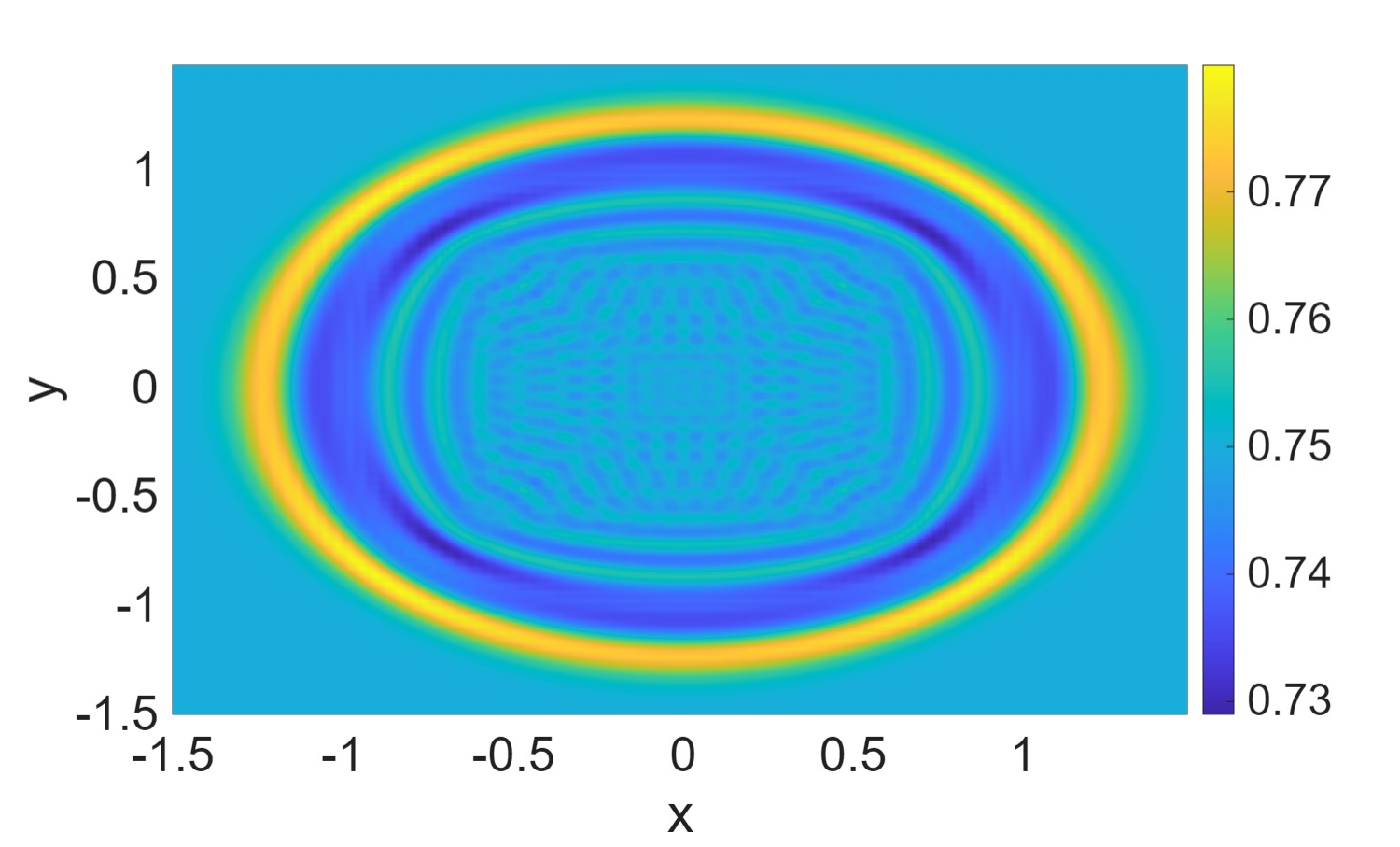}}
	\subfigure[$u(0.4,x,y)$]{\includegraphics[width=0.51\textwidth]{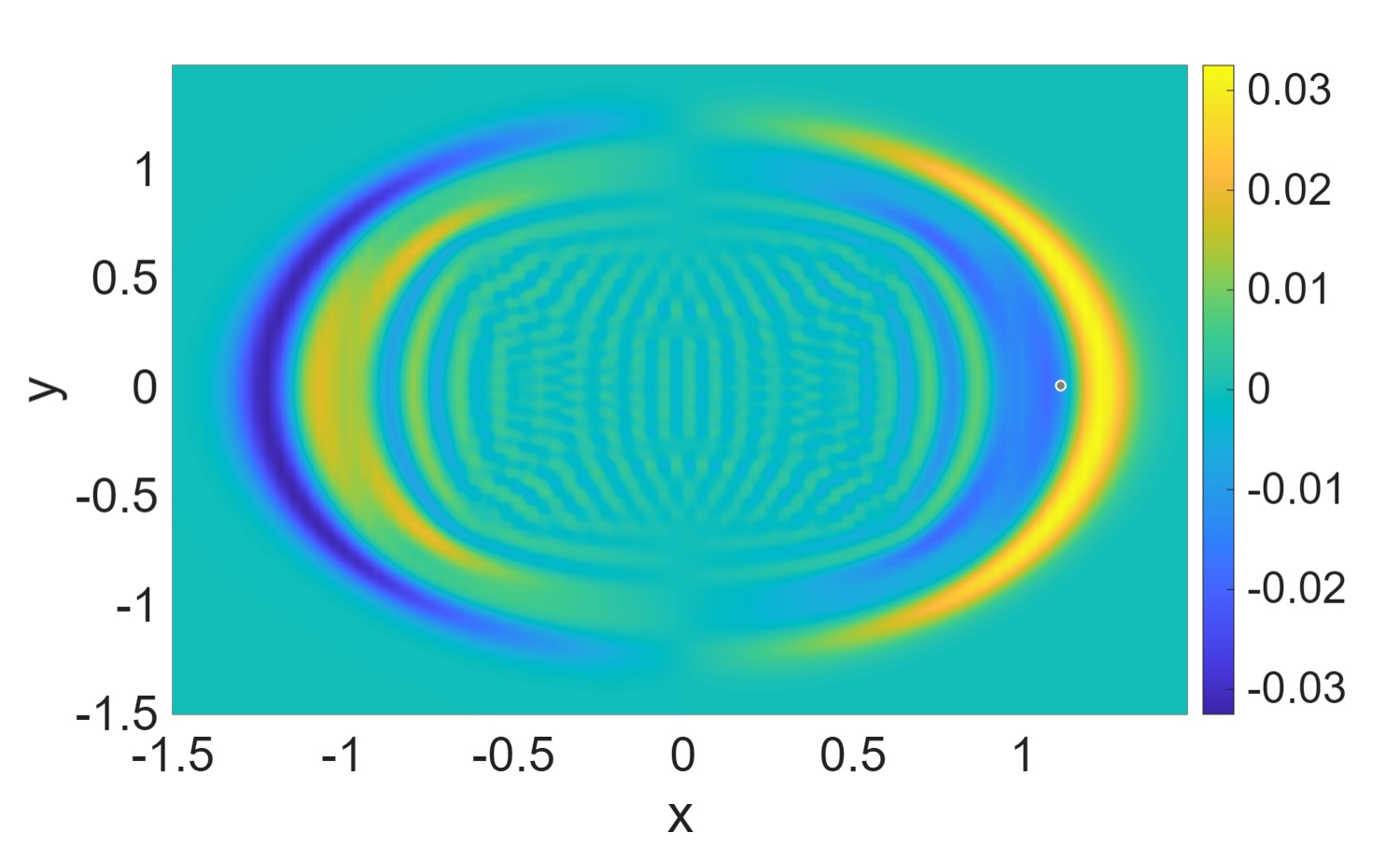}}
	\subfigure[$v(0.4,x,y)$]{\includegraphics[width=0.51\textwidth]{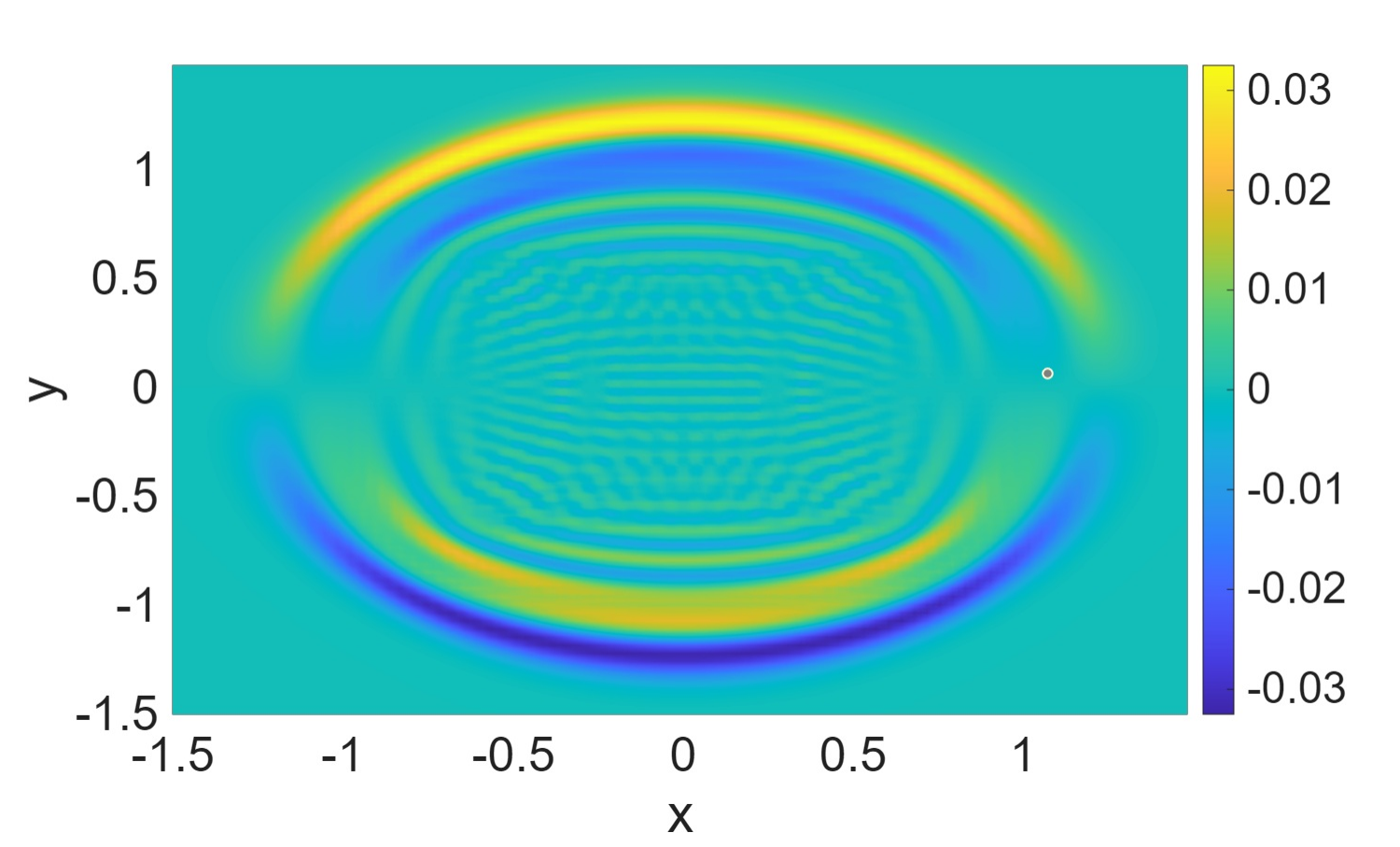}}
	\subfigure[$\omega(0.4,x,y)$]{\includegraphics[width=0.51\textwidth]{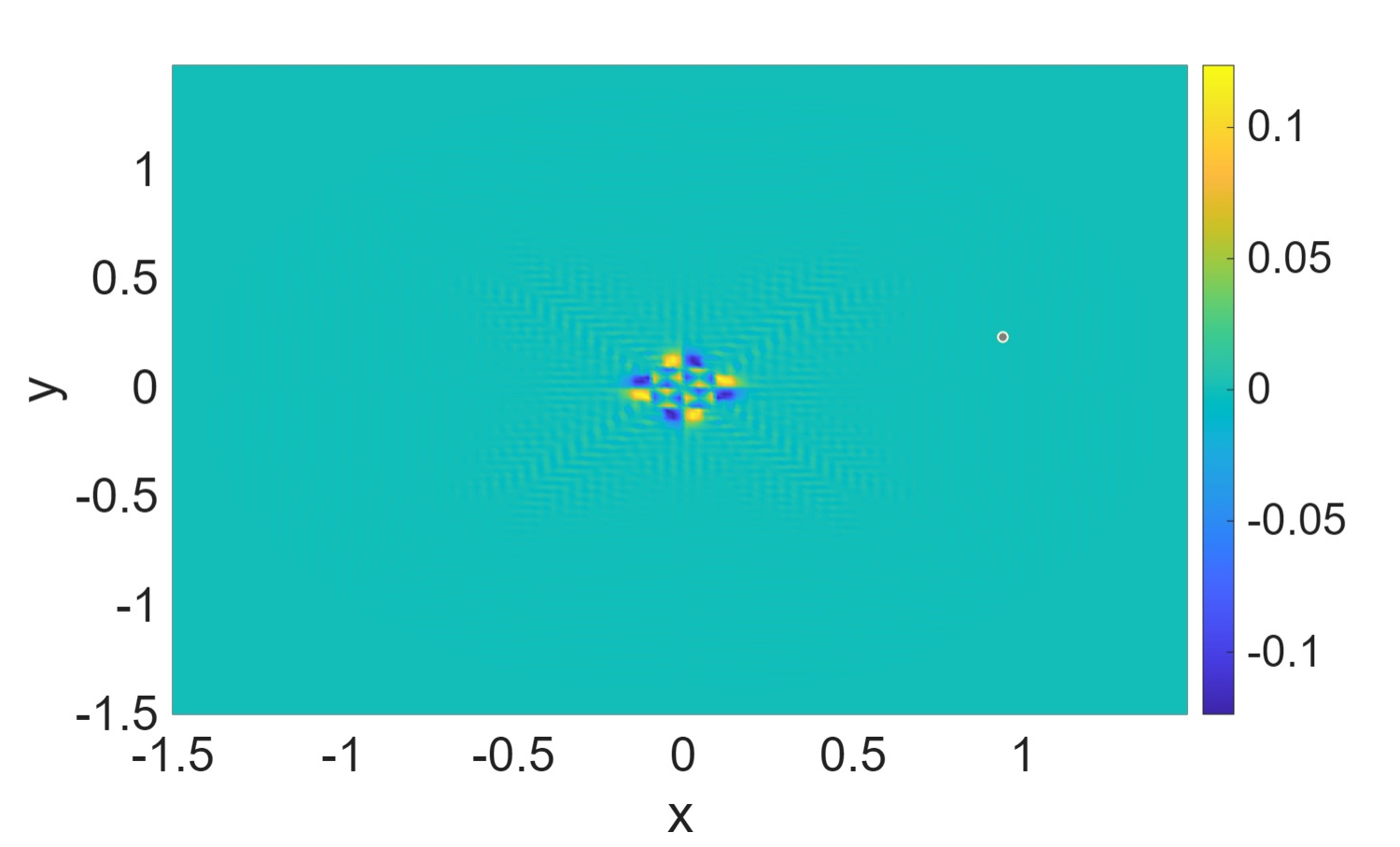}}	
	\caption{The numerical simulation of the explosive problem.}
	\label{fig5}
\end{figure}

In this section, we consider an initial overpressure that is localised in a small region. This overpressure expands and generates a propagating shock wave — in other words, an explosion. The corresponding initial condition is given by (see \cite{Lukas2022})
 $$
 \begin{aligned}
 	\rho(0, x, y) & = \begin{cases}1, & x^2+y^2 \leq R^2 \\
 		\frac{1}{10}, & \text { otherwise }\end{cases}, \\
 	u_1(0, x, y) & =0, \quad u_2(0, x, y)=0,
 \end{aligned}
 $$ 
where $(x,y) \in [-1.5,1.5]^2$ and $R = 10^{-1}$. Spatial discretization of \eqref{eq41} uses an equidistant grid with $n=160$ degrees of freedom per direction. As the spatial discretization employs standard finite differences without shock-capturing stabilization, we restrict the explosion experiment to the pre-shock regime. Beyond shock formation, the solution is no longer smooth and numerical instabilities arise. We set the kinematic viscosity to $\nu = 10^{-4}$, leading problem \eqref{eq41} becomes a strongly advection-dominated. Figure \ref{fig5} presents the numerical solutions at the final time $t = 0.4$ for the density $\rho$, the velocity components $u$ and $v$, and the vorticity $\omega = \partial_{x} v - \partial_{y} u$.

\begin{figure}[h!]
	\subfigure{\includegraphics[width=0.48\textwidth]{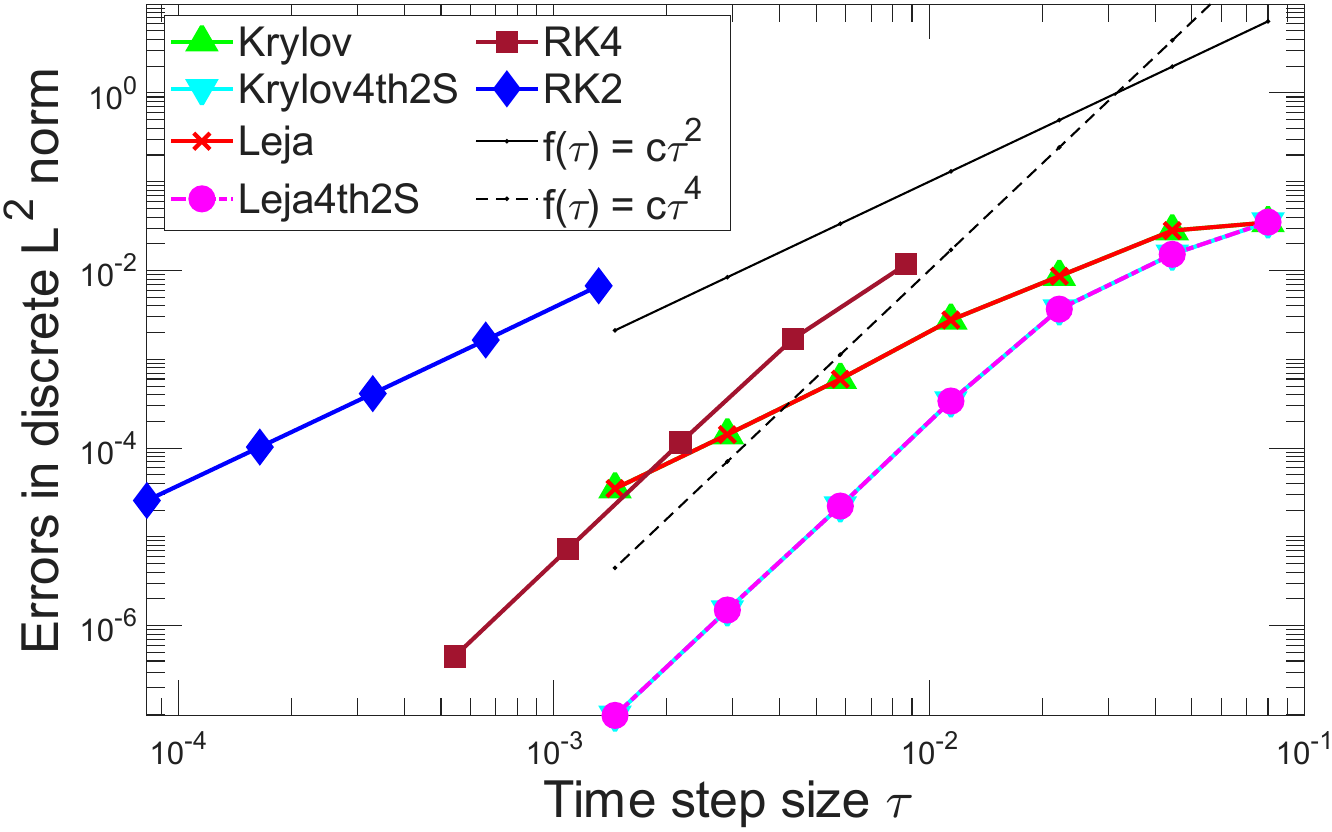}}
	\subfigure{\includegraphics[width=0.48\textwidth]{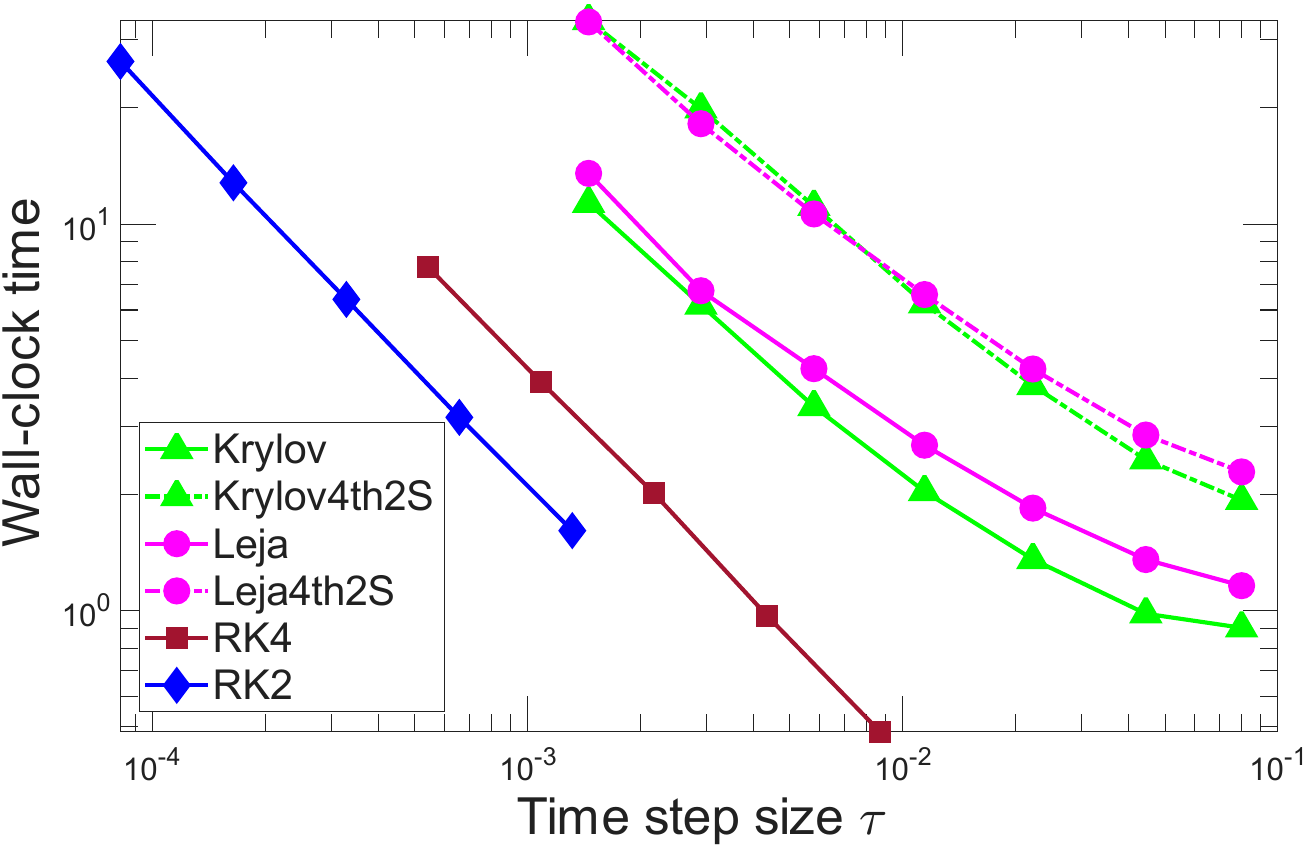}}	
	\begin{center}
		\subfigure{\includegraphics[width=0.48\textwidth]{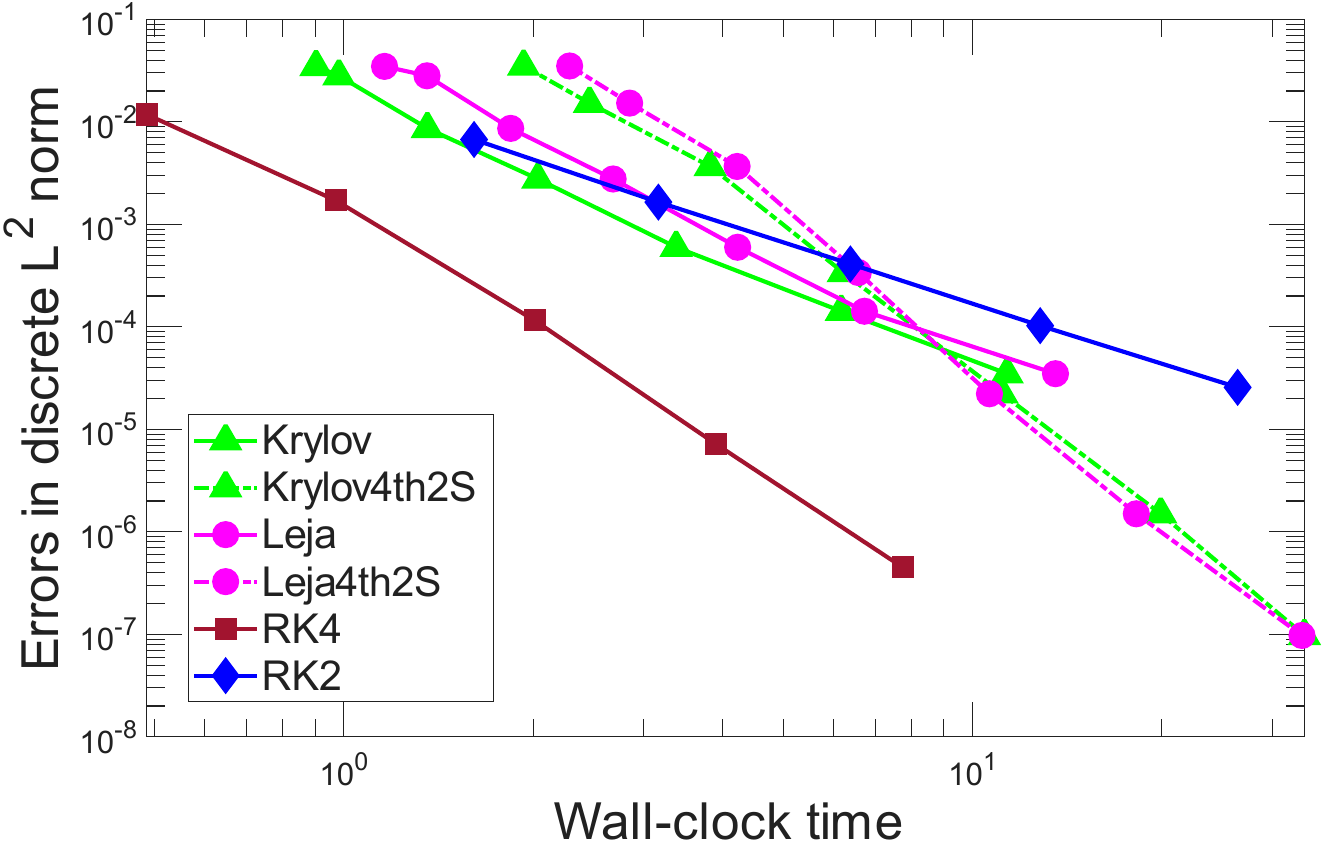}}	
	\end{center}
	\caption{The numerical results for the two-dimensional compressible isothermal Navier--Stokes problem \eqref{eq41} at time $t=0.4$. The two upper figures show the global error and the computational cost as functions of the time-step size for all methods considered. For each time-step size, the exponential Rosenbrock--Euler method yields nearly identical errors when implemented with either the Krylov or Leja approach, and the corresponding curves therefore overlap. The lower figure compares the computational cost of the different methods, where the exponential schemes are tested with time-step sizes \(\tau = \tfrac{2^{\,2-m}}{25}\) for \(1 \le m \le 7\).}	
	\label{fig6}
 
\end{figure}

Figure~\ref{fig6} presents the work–precision diagram for the methods considered in this study. In our simulation, the maximum time step size for the exponential integrators is set to $\tau = \frac{2}{25}$ and the reference solution is obtained using RK4 using a sufficiently small time step. The labels $\texttt{Krylov4th2S} $ and $\texttt{Leja4th2S} $ correspond to the  $\mathtt{exprb42}$ scheme implemented using the Krylov and Leja methods, respectively. To enhace computational efficiency, the matrix function $ \varphi_1$ is evaluated directly during the first stage of both $\mathtt{exprb42}$ the exponential Rosenbrock–Euler method. For the solution update in the $\mathtt{exprb42}$ method, the augmented matrix approach described in \cite{doi:10.1137/100788860} is employed. 

We observe that exponential integrators, whether implemented using the Krylov or Leja approach, show performance advantages similar to those observed in the linear, strongly advection-dominated setting. In particular, they allow significantly larger time-step sizes than the RK2 and RK4 schemes—by factors of approximately $61$ and $10$, respectively. However, the Krylov- and Leja-based exponential integrators achieve computational efficiency similar to that of RK2. In this test case, they are less efficient than the RK4 method.

\subsection{Shear Flow}
In this section, we consider problem \eqref{eq41} with initial data representing a shear flow. A typical configuration for this case is provided in \cite{Lukas2022,EINKEMMER2014,LIU2000577}
\begin{equation}\label{51}
\begin{aligned}
	\rho(0, x, y) & =1,	\\
	u(0, x, y) &=\begin{cases}v_0 \tanh \bigl(\frac{y-\frac{1}{4}}{d}\bigr), & y \leq \frac{1}{2}, \\
		v_0 \tanh \bigl(\frac{\frac{3}{4}-y}{d}\bigr), &  y>\frac{1}{2},\end{cases} \\
	v(0, x, y) &=\delta \sin (2 \pi x),
\end{aligned}
\end{equation}
where $(x, y) \in[0,1]^{2}$, $v_{0}=0.1$, $d=1 / 30$, and $\delta=5 \cdot 10^{-3}$. The problem under consideration describes the motion of a fluid within the domain. In our simulations, the spatial discretization employs standard finite differences on an equidistant grid with $n = 160$ degrees of freedom in each spatial direction. For a given Reynolds number $\text{Re}$, the kinematic viscosity is computed as $\nu = v_0 / \text{Re}$. In this study, we take $\text{Re} = 10^5$, which gives $\nu = 10^{-6}$, and makes the problem strongly advection-dominated.  The maximum time step size for the exponential integrators is chosen as $\tau  = 1$.

\begin{figure}[t!]
	\subfigure[$\rho(12,x,y)$]{\includegraphics[width=0.48\textwidth]{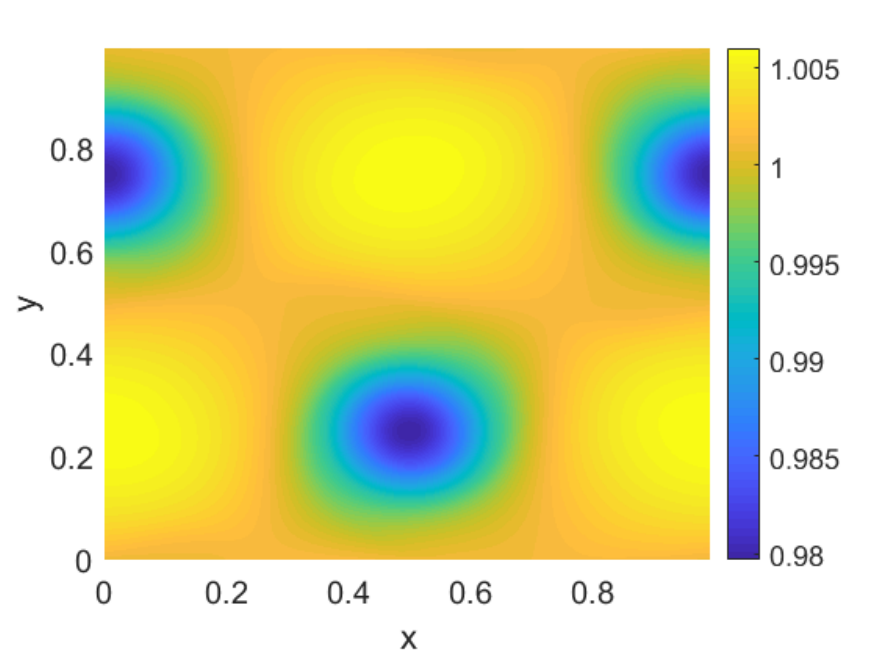}}
	\subfigure[$u(12,x,y)$]{\includegraphics[width=0.48\textwidth]{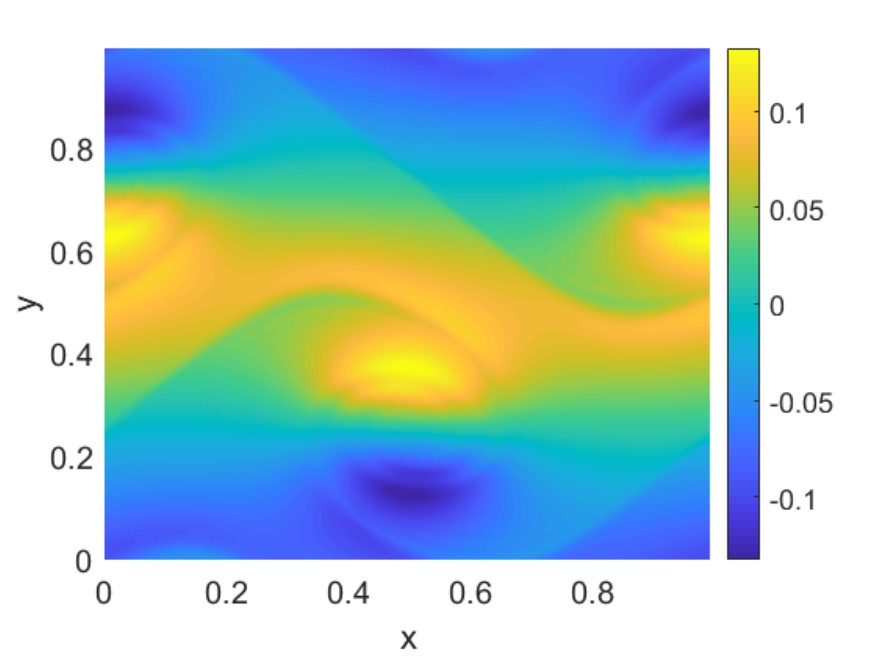}}
	\subfigure[$v(12,x,y)$]{\includegraphics[width=0.48\textwidth]{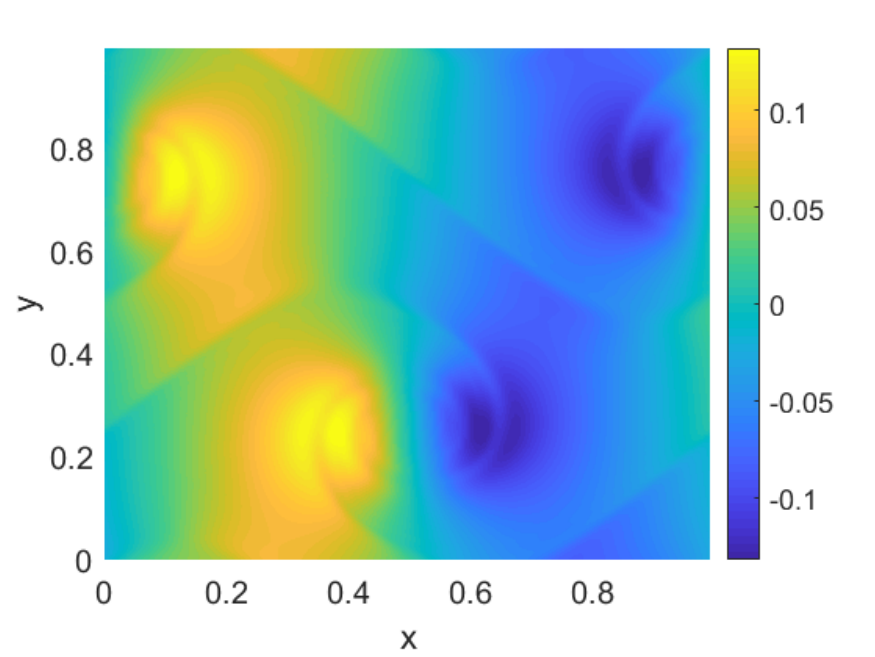}}\hspace{4.8mm}
	\subfigure[$\omega(12,x,y)$]{\includegraphics[width=0.48\textwidth]{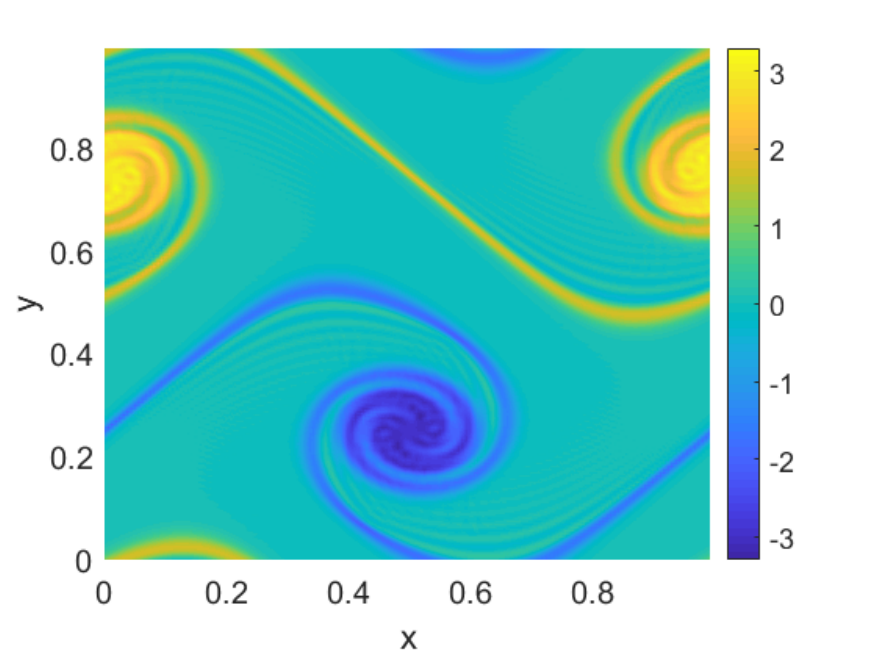}}	
	\caption{The numerical simulation of the shear flow problem.}
	\label{fig7}
\end{figure}

Figure~\ref{fig7} shows the numerical solution of density \(\rho\), velocity components \(u\) and \(v\), and vorticity  at \(t = 12\). 
Near the top and bottom boundaries, the fluid flows leftward, while it moves in the opposite direction through the central region. Over time, small perturbations in the velocity field amplify, leading to dynamic structures dominated by vortical motion.

\begin{figure}[h!]
	\subfigure{\includegraphics[width=0.51\textwidth]{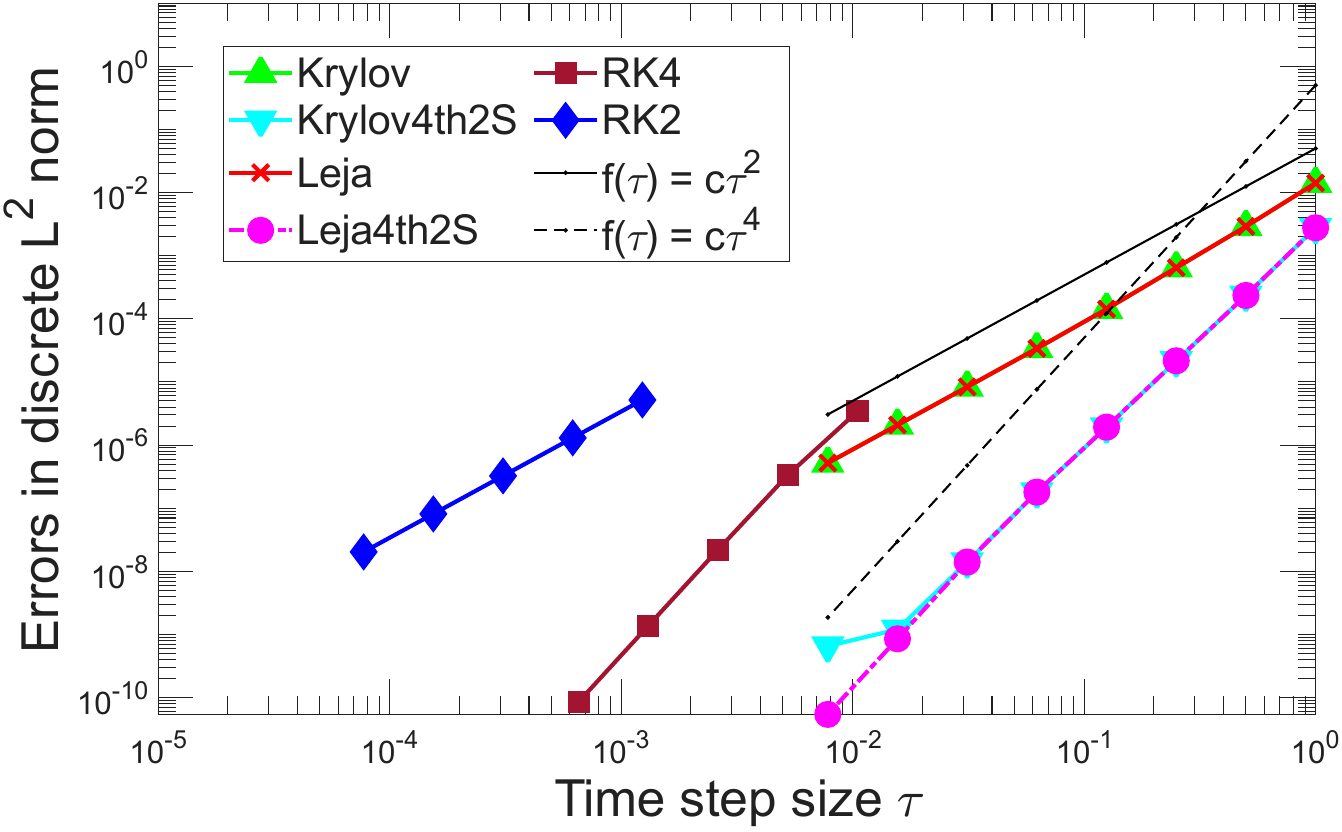}}
	\subfigure{\includegraphics[width=0.51\textwidth]{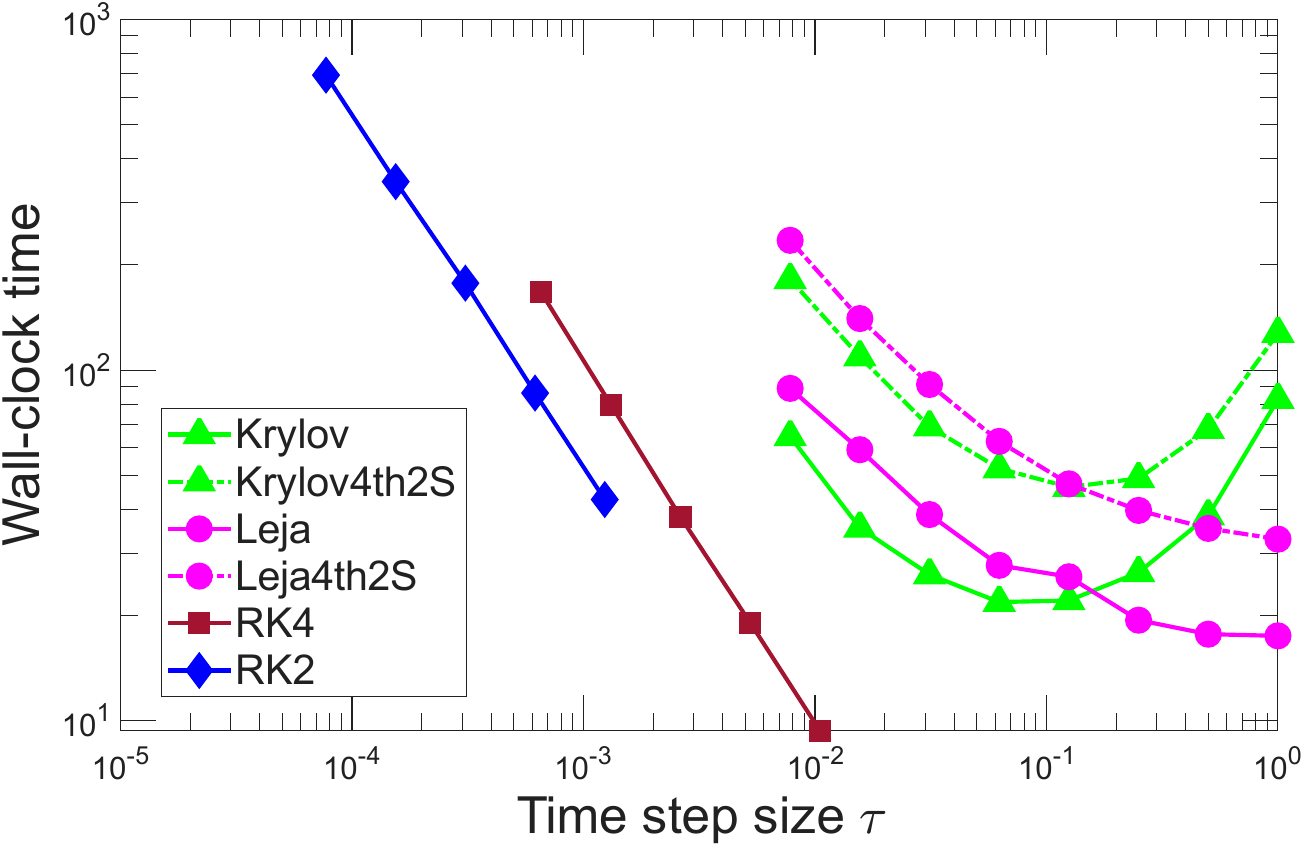}}	
	\begin{center}
		\subfigure{\includegraphics[width=0.51\textwidth]{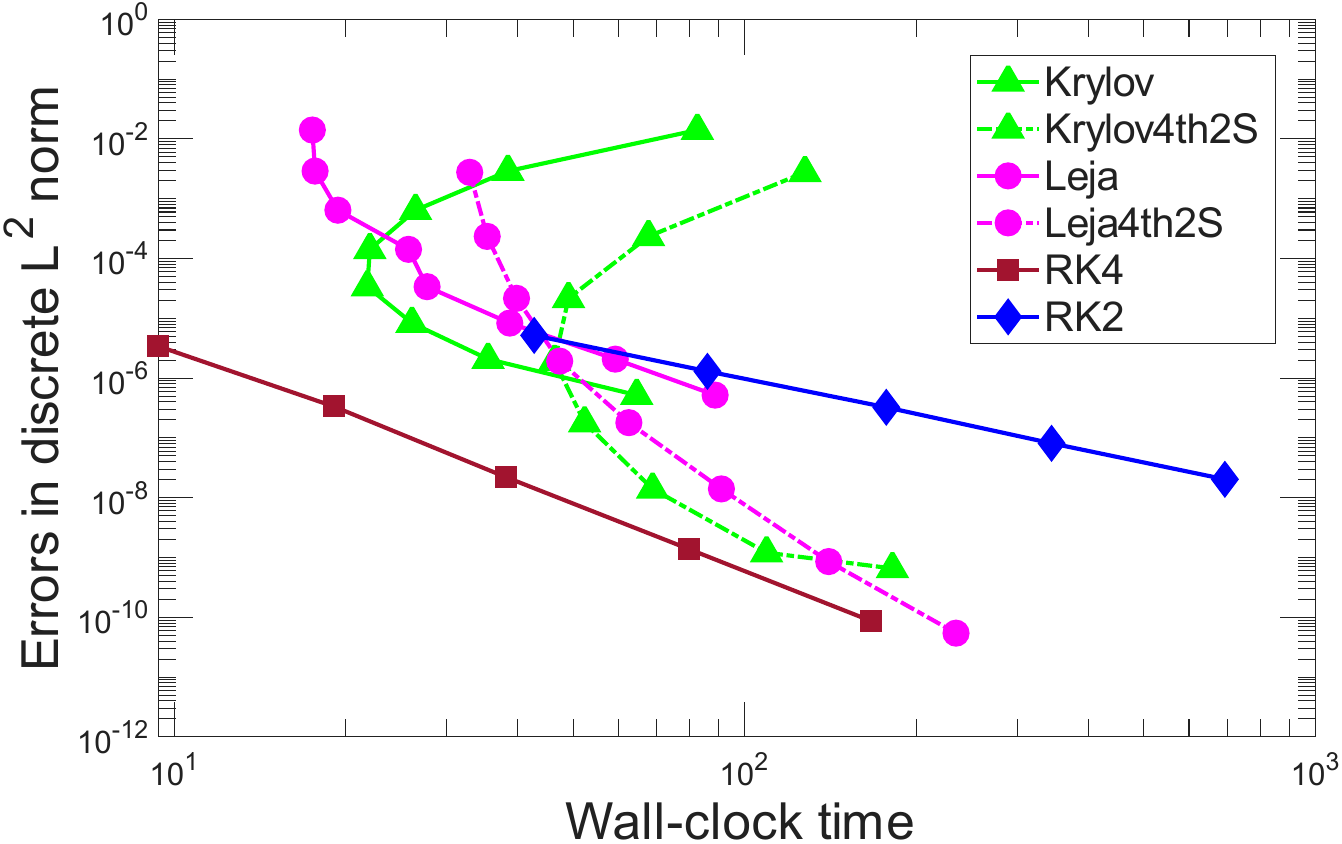}}	
	\end{center}
	\caption{The numerical results for the two-dimensional compressible isothermal Navier--Stokes problem \eqref{eq41} with Reynolds number $10^5$ at time $t=12$. The two figures at the top illustrate the global error and the computational cost as functions of the time step size for the various methods under consideration. For a given time step size, the errors of the exponential Rosenbrock--Euler method, when using the Krylov method or Leja interpolation, are similar, and the curves overlap in the figure. The bottom figure shows a comparison of the cost of the various methods, wherein the exponential methods use the time step sizes $\tau = 2^{-m}$, $0\le m\le 8$.}\label{fig3}
\end{figure}

Figure \ref{fig3}
shows the corresponding work-precision diagram. It is observed that exponential integrators allow much larger time steps than RK2 and RK4—by factors of about $810$ and $96$, respectively. The Leja-based exponential Rosenbrock–Euler method reduces the computational cost by approximately a factor of 1.7 compared to the RK2 method, while maintaining an accuracy of around \(10^{-5}\). In addition, the exponential Rosenbrock--Euler method implemented with the Krylov (with small time steps) and Leja approach  consistently outperforms RK2 across all tolerances, from moderate to stringent. 

The fourth-order \texttt{exprb42} scheme is about five times more expensive than RK4 at this accuracy. However, for a practical tolerance of \(10^{-3}\), it is only roughly twice as costly as RK4, whereas RK2 and RK4 cannot reach this regime due to stability. Moreover, \texttt{exprb42} remains competitive for smaller time steps, since halving the step size does not double the cost, unlike RK4. These results
demonstrate that exponential integrators remain competitive for all considered tolerances
and problem settings.

Although Krylov-based methods are relatively expensive for large time steps, they outperform the Leja-based method and RK2 when smaller steps are used. As in the linear case, this efficiency stems from the inner products benefiting from good cache locality.

\section{Conclusion}	\label{conclude}
We have evaluated exponential integrators for advection-dominated problems in one- and two-dimensional domains using CPU-time measurements. Across
all test problems, we observe that Leja-based methods outperform Krylov-based methods
when large time steps are employed. In contrast, Krylov methods are more efficient for small steps, outperforming Leja approaches. Overall, exponential integrators match or surpass explicit Runge--Kutta methods in accuracy and efficiency.
Exponential integrators whether Krylov- or Leja-based show clear performance advantages over explicit schemes in weakly advection-dominated problems.  In strongly advection-dominated cases, however, their performance is comparable. Their efficiency is especially notable in linear problems and in domains with mixed advection-diffusion characteristics. Large time steps are allowed in weakly advection-dominated regions, though they require more iterations per step. In linear problems, exponential integrators are particularly advantageous, achieving high accuracy without a large increase in computational cost.

\end{document}